\documentclass[12pt,a4paper]{article}
\usepackage{amsmath,amssymb,ifthen,url}
\usepackage[amsmath,thmmarks]{ntheorem}
\usepackage{braket}
\usepackage[all]{xy}
\usepackage{mathrsfs}
\SelectTips{cm}{11}

\theorembodyfont{\upshape}
\newtheorem{defn}{Definition}[section]
\newtheorem{rem}[defn]{Remark}

\theorembodyfont{\slshape}
\newtheorem{thm}[defn]{Theorem}
\newtheorem{prop}[defn]{Proposition}
\newtheorem{lem}[defn]{Lemma}
\newtheorem{cor}[defn]{Corollary}

\theoremstyle{nonumberplainwithoutbrackets}
\theorembodyfont{\upshape}
\theoremsymbol{\rule{.7em}{.7em}}
\theoremheaderfont{\itshape}
\theoremseparator{.}
\newtheorem*{prf}{Proof}

\newcommand{\Afra}{\mathfrak{A}}
\newcommand{\Bfra}{\mathfrak{B}}
\newcommand{\Kfra}{\mathfrak{K}}
\newcommand{\Pfra}{\mathfrak{P}}
\newcommand{\Qfra}{\mathfrak{Q}}
\newcommand{\ofra}{\mathfrak{o}}
\newcommand{\pfra}{\mathfrak{p}}
\newcommand{\N}{\mathbb{N}}
\newcommand{\Z}{\mathbb{Z}}
\newcommand{\Q}{\mathbb{Q}}
\newcommand{\C}{\mathbb{C}}
\newcommand{\F}{\mathbb{F}}

\DeclareMathOperator{\Hom}{Hom}
\DeclareMathOperator{\End}{End}
\DeclareMathOperator{\Aut}{Aut}
\DeclareMathOperator{\Gal}{Gal}
\DeclareMathOperator{\Ind}{Ind}
\DeclareMathOperator{\GL}{GL}
\DeclareMathOperator{\cInd}{c-Ind}
\DeclareMathOperator{\Res}{Res}
\DeclareMathOperator{\Cent}{Cent}
\DeclareMathOperator{\M}{M}
\DeclareMathOperator{\U}{\mathbf{U}}
\DeclareMathOperator{\Irr}{Irr}
\DeclareMathOperator{\I}{I}

\DeclareMathOperator{\Nrd}{Nrd}

\title{On the types for supercuspidal representations of inner forms of $\GL_N$}
\author{Yuki Yamamoto}

\begin{document}

\maketitle

\begin{abstract}
Let $F$ be a non-Archimedean local field, $A$ be a central simple $F$-algebra, and $G$ be the multiplicative group of $A$.  
It is known that for every irreducible supercuspidal representation $\pi$, there exists a $[G, \pi]_{G}$-type $(J, \lambda)$, called a (maximal) simple type.  
We will show that $[G, \pi]_{G}$-types defined over some maximal compact subgroup are unique up to $G$-conjugations under some unramifiedness assumption on a simple stratum.  

\tableofcontents

\end{abstract}

\section{Introduction}
All the representations considered in this paper are smooth and over $\C$.  
Let $F$ be a non-Archimedean local field with its integer ring $\ofra_F$, its radical $\pfra_F$ and its residue field $k_{F}=\ofra_{F}/\pfra_{F}$
Let $G$ be the multiplicative group of a central simple $F$-algebra.  
Let $\Irr(G)$ be the set of the isomorphism classes in the irreducible $G$-representations.  
Let $\mathcal{M}(G)$ be the category of representaions of $G$.  
Then there exists an equivalence relation $\sim$ on $\Irr(G)$ such that $\mathcal{M}(G)$ decomposes into indecomposable subcategories, indexed by $\mathcal{B}(G)=\Irr(G)/\sim$.  
This decomposition of category is called the Bernstein decomposition.  

\begin{defn}
Let $\mathfrak{s}$ be an element of $\mathcal{B}(G)$ and $J$ be a compact open subgroup of $G$.  
An irreducible $J$-representation $\lambda$ is called an $\mathfrak{s}$-type if the following condition holds:  
for $\pi \in \Irr(G)$, $\pi \in \mathfrak{s}$ if and only if $\pi|_{J} \supset \lambda$.  
\end{defn}

When $\pi \in \Irr(G)$ is supercuspidal, for $\pi' \in \Irr(G)$, we have $\pi \sim \pi'$ if and only if there exists an unramified character $\chi$ of $F^{\times}$ such that $\pi' \cong \pi \otimes \left( \chi \circ \Nrd_{A/F} \right)$, where $\Nrd_{A/F}$ is the reduced norm of $A$ over $F$.    
For $\pi \in \Irr(G)$ which is supercuspidal, we define
\[
	[G, \pi]_{G}=\Set{\pi' \in \Irr(G) | 
	\begin{array}{l}
	\text{$\pi' \cong \pi \otimes \left( \chi \circ \Nrd_{A/F} \right)$ for some } \\
	\text{unramified character $\chi$ of $F^{\times}$}
	\end{array}
	} \in \mathcal{B}(G).  
\]
Then $[G, \pi]_{G}$ consists of supercuspidal representations.  

Let $G = \GL_{N}(F)$ and $K$ be a maximal compact subgroup of $G$.  
For $N=2$, Henniart proved in the Appendix to \cite{BM} that if $\pi$ is an irreducible supercuspidal representation of $G$, then there exists a unique representation $\rho$ of $K$ such that a pair $(K, \rho)$ is $[G, \pi]_{G}$-type.  
Henniart also determined the number of isomorphism classes of irreducible representation $\rho$ of $K$ such that $(K, \rho)$ is an $\mathfrak{s}$-type for every $\mathfrak{s} \in \mathcal{B}(G)$.  
Paskunas extended the Henniart's result on supercuspidal representations for general $N$ in his paper \cite{Pas}, using the concept of simple types for $\GL_{N}(F)$, established in \cite{BK1}.  
However, for general $G$, it is not necessarily so.  
When $D \neq F$, a maximal compact subgroup $K$ is properly contained in the normalizer $N_{G}(K)$ of $K$ in $G$.  
If $(K, \rho)$ is $[G, \pi]_{G}$-type and $k \in N_{G}(K)$, then $(K, {^k}\rho)$ is also $[G, \pi]_{G}$-type, where ${^k}\rho$ is the $K$-representation obtained by twisting $\rho$ by $k$.  
For some $k \in N_{G}(K) \backslash K$, it may happen that ${^k}\rho$ is not isomorphic to $\rho$.  

Then, we define $\mathscr{T}(\pi)$ to be the set of $(K, \tau)$ such that {$K$ is a maximal compact subgroup of $G$, $\tau$ is a $K$-representation and $(K, \tau)$ is a $[G, \pi]_{G}$-type.  
Every $G$-conjugacy class of $\mathscr{T}(\pi)$ is called a $[G, \pi]_{G}$-archetype.  

The main result of this paper is following:

\begin{thm}
Let $G=GL_{m}(D)$ and $\pi$ be an irreducible supercuspidal representation of $G$.  
If $\pi$ contains a simple type $(J, \lambda)$ attached to the simple stratum $[\Afra, n, 0, \beta]$ such that $F[\beta]$ is an unramified extension of $F$, then there exists a unique $[G, \pi]_{G}$-archetype.  
\end{thm}

Our approach is similar to that in \cite{Pas}.  
However, the same approach as in \cite{Pas} cannot be applied in our case.  

In the proof of the main theorem, we consider depth-zero supercuspidal representations of $G$.  
One of the key point of our proof is that for a given depth-zero supercuspidal repsentation $\sigma$, there exists some representation $\sigma'$ which is 'similar to' $\sigma$ in some sense.  
For our purpose, we cannot take representations which intertwine with $\sigma$ as the candidate for $\sigma'$.  
If $D=F$, such a representation is nothing but $\sigma$.  
However, when the dimension of $D$ over $F$ increases, the number of such representations also increases.  
Therefore, to show the existence of the representation $\sigma'$ with desired condition, we need to examine depth-zero supercuspidal representations more closely than \cite{Pas}.  

In fact, we will construct some irreducible supercuspidal representation $\pi$ of $G$ such that a $[G, \pi]_G$-archetype is not unique in some case.  
Then, this representation is a counterexample for the conjecture on the unicity of archetypes in \cite[Conjecture 4.4]{La}.  
This conjecture holds for the split case $G=\GL_N(F)$ and the multiplicative group $G=D^{\times}$ of a central division algebra over $F$ (for $G=D^{\times}$, see, for example, \cite{Do}).
It is interesting that for the multiplicative group $G=\GL_m(D)$ of a general central simple algebra over $F$, which is similar to these groups, the conjecture does not hold any longer.  

In \cite{LN}, the unicity of types is discussed for tame, toral and positive regular supercuspidal representations of a simply connected and semisimple group.  
Latham--Nevins showed that if we assume some unramifiedness on the representation, then the conjecture in \cite{La} holds.  
For $G=\GL_{m}(D)$, a simple stratum $[\Afra, n, 0, \beta]$ for a toral supercuspidal representation is considered to satisfy the condition that $F[\beta]$ is a maximal subfield of $A$.  
From this point of view, both of the unramified assumptions seem to be linked to each other.  
We note that the result in this paper covers some nontoral representations of $G$.  

We sketch the outline of this paper.  
In \S\ref{prelim}, we introduce the consepts which we will need later.  
In \S\ref{decomp}, using the Mackey decomposition we consider the restriction of an irreducible supercuspidal representation $\pi$ to a maximal compact open subgroup $K$ of $G$.  
We need to consider the double cosets of $K \backslash G / \Kfra(\Afra)$, where $\Kfra(\Afra)$ is the normalizer of some hereditary $\ofra_F$-order $\Afra$ in $A$, defined in \S\ref{prelim}.  
In \S\ref{nice_repre}, for every $g' \in G \setminus K\Kfra(\Afra)$ we take a nice representative $g \in Kg'\Kfra(\Afra)$.  
In \S\ref{propertyA}, we consider irreducible $K$-representations contained in the representations corresponding to a pair $\left(g, Kg'\Kfra(\Afra) \right)$ in the Mackey decomposition of $\pi|_{K}$ when the pair $\left(g, Kg'\Kfra(\Afra) \right)$ has `Property A'.  
In \S\ref{mainthm}, the main theorem of this paper is proved.  

\bigbreak

\noindent{\bfseries Acknowledgment}\quad
I would like to thank my supervisor Naoki Imai for his tremendous support, constant encouragement and helpful advice.  
He also checked the draft of this paper and pointed out mistakes.  
I am supported by the FMSP program at Graduate School of Mathematical Sciences, the Univercity of Tokyo.  

\bigbreak

\section{Notation and Preliminary}
\label{prelim}

Let $F$ be a non-Archimedean local field. 
Let $\ofra _F$ be its integers, $\pfra _F$ be the prime ideal of $\ofra _F$ and $k_F=\ofra_F / \pfra_F$ be the residue field of  $F$ with the cardinality $q_F$.  
If $D$ is a finite extension field of $F$ or a division $F$-algebra, we use the same notation $\ofra_D, \pfra_D, k_D, q_D$.  

\subsection{Intertwining}
If $H$ is a subgroup of $G$ and $g$ is an element of $G$, then we denote ${^{g}}H = gHg^{-1}$.  
Moreover, if $\sigma$ is a representations of $H$, we define the representation $^{g}\sigma$ of $^{g}H$ as follows:
\[
	^{g}\sigma(x) = \sigma(g^{-1}xg), \quad x \in {^{g}H}.  
\]
Let $H_1, H_2, G'$ be subgroups of $G$, $H$ be a subgroups of $H_1 \cap H_2$.  
Let for $i=1,2$, $\sigma_i$ be a representation of $H_i$ such that $\Hom _{H}(\sigma_1, \sigma_2)$ is finite-dimensional.  
We define
\begin{eqnarray*}
	\left< \sigma_1, \sigma_2 \right>_{H} &=& \dim \Hom_{H}(\sigma_1, \sigma_2), \\
	I_{g}(\sigma_1, \sigma_2) &=& \Hom_{H_1 \cap {^g}H_2}(\sigma_1, {^g}\sigma_2) \text{ for } g \in G, \\
	I_{G'}(\sigma_1, \sigma_2) &=& \Set{g \in G'| I_{g}(\sigma_1, \sigma_2) \neq 0 }.  
\end{eqnarray*}
We say $\sigma_1$ and $\sigma_2$ intertwine in $G$ if $I_{G}(\sigma_1, \sigma_2) \neq \emptyset$.  
When $H_1$ and $H_2$ are compact, $\sigma_1$ and $\sigma_2$ intertwine in $G$ if and only if $\sigma_2$ and $\sigma_1$ intertwine.  

\subsection{Lattices, hereditary orders}

We recall difinitions on lattices and hereditary orders from \cite[\S 1]{S3}.  

Let $A$ be a central simple $F$-algebra and $G$ be the multiplicative group of $A$.    
Let $V$ be a simple left $A$-module.  
Then, $\End_{A}(V)$ is a central division $F$-algebra.  
Let $D$ denote the opposite of $\End_{A}(V)$.  
Then, $V$ becomes a right $D$-vector space with $\dim_{D}V=m$ and the natural $F$-algebra isomorphism $A \cong \End_{D}(V)$ exists.  
If we fix a right $D$-isomorphism $V \cong D^m$, this isomorphism induces $A \cong \M_{m}(D)$ and $G \cong \GL_{m}(D)$.  

Let $E/F$ be a subfield in $A$.  
Then $V$ is equipped with an $E$-vector space structure which is compatible with the $D$-action, whence $V$ becomes a right $E \otimes_F D$-module.  
Let $W$ be a simple right $E \otimes _F D$-module.  
Then, there exists $m' \in \N$ and a right $E \otimes _F D$-module isomorphism $V \cong W ^{\oplus m'}$.  
Put $A(E)=\End_D (W)$.  
By the above isomorphism, we have $F$-algebra isomorphisms $\End_{D}(V) \cong \End_{D} (W ^{\oplus m'}) \cong \M_{m'} \left( \End_{D}(W) \right) = \M_{m'} \left( A(E) \right)$.  
The centralizer $\Cent_{A(E)}(E)$ of $E$ in $A(E)$ is the set of $D$-morphisms which are also $E$-linear maps, that is, $E \otimes _F D$-linear maps.  
Therefore, $\Cent_{A(E)}(E)=\End_{E \otimes _F D}(W)$ is a division $E$-algebra, whose opposite algebra is denoted by $D'$.  
Let $B$ be the centralizer of $E$ in $A$.  
Under the isomorphism $A \cong \M_{m'} \left( A(E) \right)$ , $E$ is diagonally embedded in $\M_{m'} \left( A(E) \right)$.  
Therefore, the isomorphism maps $B$ to $\M_{m'}(D')$.  

An $\ofra _D$-submodule $\Lambda$ in $V$ is called an $\ofra_D$-lattice in $V$ if it is a compact open submodule. 

\begin{defn}
Let $\Lambda_{i}$ be an $\ofra_D$-lattice in $V$ for every $i \in \Z$.  
We say that $\Lambda = \left( \Lambda_{i} \right) _{i \in \Z}$ is an $\ofra_D$-sequence if 
\begin{enumerate}
\item for all $i \in \Z$, we have $\Lambda_{i} \supset \Lambda_{i+1}$, 
\item there exists $e = e(\Lambda) \in \N$ such that $\Lambda_{i+e} = \Lambda_{i} \pfra_D$ for all $i \in \Z$.  
\end{enumerate}
The number $e=e(\Lambda)$ is called the period of $\Lambda$.  
An $\ofra_D$-sequence $\Lambda = \left( \Lambda_{i} \right)$ in $V$ is called an $\ofra_D$-chain if $\Lambda_{i} \supsetneq \Lambda_{i+1}$ for every $i$.  
\end{defn}

Let $\Afra$ be an $\ofra_F$-order in $A$.  
Then, $\Afra$ is hereditary if every left and right ideal in $\Afra$ is $\Afra$-projective.  

Let $\Lambda = \left( \Lambda_{i} \right)$ be an $\ofra_D$-sequence in $V$.  
Put $\Pfra_{i}(\Lambda)=\Set{ x \in A | x\Lambda_{j} \subset \Lambda_{i+j}, j \in \Z }$.  
Then, $\Afra = \Pfra_{0}(\Lambda)$ is a hereditary $\ofra_F$-order in $A$.  
The radical of $\Afra$ is $\Pfra(\Afra)=\Pfra_{1}(\Lambda)$.  
For every hereditary $\ofra_F$-order $\Afra$, there exists an $\ofra_D$-chain $\Lambda$ such that $\Afra = \Afra(\Lambda)$.  
If $[\Lambda_i : \Lambda_{i+1}]$ is constant for any $i$, then $\Afra=\Afra(\Lambda)$ is called principal.  

Let $\Lambda=\left( \Lambda_{i} \right)$ be an $\ofra_D$-chain in $V$.  
Let $\Kfra(\Lambda)$ be the set of $g \in A$ with the condition that there exists $n \in \N$ such that $g(\Lambda_{i})=\Lambda_{n+i}$ for all $i$.  
For the hereditary $\ofra_F$-order $\Afra = \Afra(\Lambda)$, we put $\Kfra(\Afra)=\Set{ g \in G | g\Afra g^{-1}=\Afra }$.  
Then $\Kfra(\Afra)$ is equal to $\Kfra(\Lambda)$ and the group homomorphism $v_{\Afra} \colon \Kfra(\Afra) \to \Z$ is defined as $v_{\Afra}(g)=n$ for $g \in \Kfra(\Afra)$ with $g\Afra=\Pfra(\Afra)^{n}$.  
Then, the short exact sequence
\[
	1 \to \U (\Afra) \to \Kfra (\Afra ) \to \Z
\]
exists, where $\U (\Afra) = \U^{0}(\Afra) = \Afra ^{\times}$ is the unique maximal compact open subgroup of $\Kfra (\Afra )$.  
For $n \in \N_{>0}$, put $\U^{n}(\Afra) = 1+\Pfra(\Afra)^{n}$.  

A hereditary $\ofra_F$-order $\Afra$ in $A$ is $E$-pure if $E^{\times} \subset \Kfra(\Afra)$.  

\begin{thm}[{{\cite[Theorem 1.3]{Bro}}}]
Let $\Afra$ be an $E$-pure hereditary $\ofra_F$-order in $A$ with its radical $\Pfra$.  
Then, $\Bfra = \Afra \cap B$ is a hereditary $\ofra_E$-order in B with its radical $\Qfra=\Pfra \cap B$.  
\end{thm}

If $A(E)$ is as above, there is a unique $E$-pure hereditary $\ofra_{F}$-order $\Afra(E)$ with its radical $\Pfra(E)$; see \cite[\S 1.5.2]{S3}.  
The hereditary $\ofra_{F}$-order $\Afra (E)$ is principal.  
The equations $\Afra(E) \cap D'=\ofra_{D'}$ and $\Pfra(E) \cap D' =\pfra_{D'}$ hold.  

\subsection{Proper hereditary orders}
Let $A, E, B$ be as above.  
If $V$ is a simple left $A$-module and $W$ is a simple right $E \otimes _F D$-module, then an $E \otimes _F D$-module isomorphism $V \cong W^{\oplus m'}$ induces the $F$-algebra isomorphism $A \cong \M_{m'} \left( A(E) \right)$.  
\begin{defn}[{{\cite[\S 4.3.1]{S3}}}]
A hereditary $\ofra_F$-order $\Afra$ in $A$ is proper if there exists an $E \otimes _F D$-module isomorphism $V \cong W^{\oplus m'}$ and integers $r, s$ such that $rs=m'$ holds and 
\[
	\Afra =
	\begin{pmatrix}
	\M_s \left( \Afra (E) \right) & \cdots & \M_s \left( \Afra (E) \right) \\
	\vdots & \ddots & \vdots \\
	\M_s \left( \Pfra (E) \right) & \cdots & \M_s \left( \Afra (E) \right) \\
	\end{pmatrix}
\]
under the identifaction between $A$ and $\M_{m'} \left( A(E) \right) = \M_r \left( \M_s \left( A(E) \right) \right)$ given by the isomorphism $V \cong W^{\oplus m'}$.  
\end{defn}

If $\Afra$ is as above, since $\Afra (E)$ is principal, $\Afra$ is also principal and 
\[
	\Bfra = \Afra \cap B =
	\begin{pmatrix}
	\M_s \left( \ofra_{D'} \right) & \cdots & \M_s \left( \ofra_{D'} \right) \\
	\vdots & \ddots & \vdots \\
	\M_s \left( \pfra_{D'} \right) & \cdots & \M_s \left(  \ofra_{D'} \right) \\
	\end{pmatrix}
\]
holds, whence $\Bfra$ is a principal hereditary $\ofra_{E}$-order in $B$ with its period $r$.  

\begin{rem}
\label{maxproper}
\begin{enumerate}
\item If $\Afra$ is proper and $\Bfra$ is maximal, then $r=1$ and we have an isomorphism $\Afra \cong \M_s \left( \Afra (E) \right)$.  
\item If $\Bfra$ is a principal hereditary $\ofra_{E}$-order in $B$, then there exists a proper hereditary $\ofra_F$-order $\Afra$ in $A$ such that $\Afra \cap B = \Bfra$, see \cite[Remarque 4.8]{S3}.  
It is proved in \cite[Lemme 1.6]{S2} that if $\Bfra$ is maximal, then there is a unique hereditary $\ofra_F$-order $\Afra$ in $A$ such that $\Afra \cap B = \Bfra$, whence $\Afra$ is proper.  
\end{enumerate}
\end{rem}

\subsection{Strata, simple characters, simple types}

\begin{defn}[{{\cite[\S 2.1]{S3}}}]
A 4-tuple $\left[ \Afra , n, r, \beta \right]$ is called a stratum of $A$ if the following conditions hold:
\begin{itemize}
\item $\Afra$ is a hereditary $\ofra _F$-order in A.  
\item $n, r$ are integers with $n > r$.  
\item $\beta$ is an element of $\Pfra(\Afra)^{-n}$.   
\end{itemize}
A stratum $\left[ \Afra , n, r, \beta \right]$ is pure if 
\begin{enumerate}
\item $E=F[\beta]$ is a subfield, 
\item $\Afra$ is $E$-pure, 
\item $v_{\Afra}(\beta) = -n$.  
\end{enumerate} 
A pure stratum $\left[ \Afra , n, r, \beta \right]$ is simple if  $r < -k_{0}(\beta, \Afra)$, where for $\beta \in A$ the element  $k_{0}(\beta, \Afra) \in \Z \cup \{-\infty\}$ is defined as in \cite[\S 2.1]{S3} and satisfies $v_{\Afra}(\beta) \leq k_{0}(\beta, \Afra)$.  
\end{defn}

If $\left[ \Afra , n, r, \beta \right]$ is a simple stratum of $A$, we can define subgroups $J(\beta, \Afra)$ and $H(\beta, \Afra)$ of $\U(\Afra)$ as in \cite[\S 3]{S1}.  
For $i \in \N$, put $J^i(\beta, \Afra) = J(\beta, \Afra) \cap \U^i(\Afra)$ and $H^i(\beta, \Afra)=H(\beta, \Afra) \cap \U^i(\Afra)$.  
We have an equation $J(\beta, \Afra) = \U(\Bfra) J^1(\beta, \Afra)$, where $\Bfra = \Afra \cap B$.  

We can also define the set of `simple characters' $\mathscr{C}(\beta, m, \Afra)$, consisting of characters of $H^1(\beta, \Afra)$ satisfying some condition, initially defined in \cite[D\'efinition 3.45]{S1}.  
We put $\mathscr{C}(\beta, \Afra) = \mathscr{C}(\beta, 0, \Afra)$.  

We recall properties on $\mathscr{C}(\beta, \Afra)$ from \cite{S2}.  
For every $\theta \in \mathscr{C}(\beta, \Afra)$, there exists a unique irreducible representation $\eta _{\theta}$ of $J^1(\beta, \Afra)$ such that $\eta _{\theta}|_{H^1(\beta, \Afra)}$ is a direct sum of $\theta$, called the Heisenberg representation of $\theta$.  
Moreover, for every $\eta_{\theta}$ as above, there exist a number of extensions $\kappa$ to $J(\beta, \Afra)$ such that $B^{\times} \subset I_{G}(\kappa, \kappa)$, called the $\beta$-extensions of $\eta_{\theta}$.  
If $\kappa$ is as above, then $I_{G}(\eta_\theta, \eta_\theta)=I_{G}(\kappa, \kappa)=JB^{\times}J$ and for $g \in JB^{\times}J$,  we have $\dim I_{g}(\eta_\theta, \eta_\theta)=\dim I_{g}(\kappa, \kappa)=1$.  

\begin{defn}[{{\cite[\S 2.5 and \S 4.1]{S3}}}]
\label{def_of_simple_type}
Let $J$ be a compact open subgroup of $G$ and $\lambda$ be an irreducible $J$-representation.  
A pair $(J, \lambda)$ is called a simple type of level 0 if
\begin{enumerate}
\item for some principal hereditary $\ofra_F$-order $\Afra$, we have $J=\U (\Afra)$, 
\item there exists an irreducible cuspidal representation $\sigma$ of $ \GL_{s}(k_{D})$ such that $\lambda$ is the lift of the irreducible representation $\sigma ^{\otimes r}$ of $\U(\Afra)/\U^1(\Afra) \cong \GL_{s}(k_{D})^{r}$ to $\U(\Afra)$.  
\end{enumerate}
A pair $(J, \lambda)$ is called a simple type of level $>0$ if
\begin{enumerate}
\item there exists a simple stratum $[ \Afra, n, 0, \beta]$ with $n>0$ such that $J=J(\beta, \Afra)$, 
\item there exist irreducible representations $\kappa, \sigma$ of $J$ with $\lambda = \kappa \otimes \sigma$ such that
\begin{itemize}
\item $\kappa$ is a $\beta$-extension of the Heisenberg representation $\eta_{\theta}$ for some $\theta \in \mathscr{C}(\beta, \Afra)$, and
\item $\sigma$ is a $J(\beta, \Afra)$-representation trivial on $J^1(\beta, \Afra)$, and when $J(\beta, \Afra)$-representations which are trivial on $J^1(\beta, \Afra)$ are regarded as $\U (\Bfra)$-representations under the group isomorphism 
\[
J(\beta, \Afra)/J^1(\beta, \Afra) \cong \U (\Bfra)/ \U ^1(\Bfra), 
\]
$(\U (\Bfra), \sigma)$ is a simple type of level 0, where $B = \Cent_{A}(F[\beta]), \Bfra = \Afra \cap B$.  
\end{itemize}
\end{enumerate}
A pair $(J, \lambda)$ is a simple type if it is a simple type of level 0 or level $>0$.  

\end{defn}

\subsection{Maximal simple types}
\label{max_simple_types}

\begin{defn}[{{\cite[5.1]{S3}}}]
Let $(J, \lambda)$ be a simple type.  
The simple type $(J, \lambda)$ is called maximal if one of the following conditions holds: 
\begin{enumerate}
\item The simple type $(J, \lambda)$ is of level 0 and $J=\U(\Afra)$ for some maximal $\ofra_{F}$-order $\Afra$ in $A$,  
\item The simple type $(J, \lambda)$ is of level $>0$, attached to the simple stratum $[\Afra, n, 0, \beta]$ and $\Bfra = \Afra \cap \Cent_{A}(F[\beta])=\Afra \cap B$ is a maximal $\ofra_{F[\beta]}$-order in $B$.  
\end{enumerate}
\end{defn}

\begin{thm}[{{\cite[Theorem 5.5(ii)]{GSZ}}} and {{\cite[Th\'eoreme 5.21]{SS}}}]
Let $(J, \lambda)$ be a simple type.  
The simple type $(J, \lambda)$ is maximal if and only if there exists an irreducible supercuspidal representation $\pi$ of $G$ such that $\lambda$ is contained in $\Res_{J}^{G}\pi$.  
If $(J, \lambda)$ is comtained in some irreducible supercuspidal representation $\pi$ of $G$, then $(J, \lambda)$ is a $[G, \pi]_{G}$-type.  
If $\pi$ is an irreducible supercuspidal representation, then there exists a simple type $(J, \lambda)$ such that $\lambda \subset \pi|_{J}$.  
\end{thm}

For a maximal simple type $(J, \lambda)$, we define a subgroup $\tilde{J}$ of G, containing and normalizing $F^{\times}J$.  

Let $\Afra$ be a maximal $\ofra_{F}$-order in $A$ in this paragraph.  
We fix an isomorphism $A \cong \M_{m}(D)$ such that under the identification by the isomorphism, $\Afra = \M_{m}(\ofra_D)
$.  
This isomorphism induces $\U(\Afra)/\U^{1}(\Afra) \cong \GL_{m}(k_D)$.  
Let $\sigma$ be a representation of $\U(\Afra)$, which is trivial on $\U^{1}(\Afra)$.  
Then $\sigma$ is the lift of a representation $\tau$ of $\GL_{m}(k_D)$ to $\U(\Afra)$.  
For $\gamma \in \Gal(k_D/k_F)$, we define
\[
	\left( \gamma \cdot \sigma \right) (g) = \tau \left( \left(\gamma^{-1} (g_{i,j} \text{ mod } p_{D} ) \right)_{i,j}  \right), \quad g=(g_{ij})_{i,j} \in \U(\Afra).  
\]

When $(J, \lambda)$ is of level 0, let $\Afra$ be a hereditary $\ofra_F$-order in $A$ such that $J=\U(\Afra)$, which is a maximal compact subgroup in $G$, and fix an isomorphism $A \cong \M_{m}(D)$ with $\Afra = \M_{m}(\ofra_D)$.  
We put
\[
	l_0=\min \Set{n' \in \N_{>0} | \gamma ^{n'} \cdot \sigma \cong \sigma \text{ for any } \gamma \in \Gal(k_D/k_F)}.  
\]
We also fix a uniformizer $\varpi_{D}$ of $D$.  
Let $y$ be the diagonal matrix of $\M_m(D)$ with every diagonal entry $\varpi_{D}$.  
Then $\tilde{J} = I_G(\sigma, \sigma)$ is a subgroup of $G$ generated by $\U(\Afra)$ and $y^{l_0}$.  

When $(J, \lambda)$ is of level $>0$, let $[\Afra, n, 0, \beta]$ be a simple stratum giving $(J, \lambda)$.  
Then there exist $\kappa, \sigma$ as Definition \ref{def_of_simple_type}.  
Since $(J, \lambda)$ is maximal, $\Bfra$ is a maximal hereditary $\ofra_E$-order in $B$ and $(\U(\Bfra), \sigma)$ is a maximal simple type of $B$ of level 0.  
Therefore we can choose $y \in B^{\times}$ and $l_0$ as well as the case of simple type of level 0.  
Then $\tilde{J} = I_G(\lambda, \lambda)$ is a subgroup of $G$, generated by $J, F^{\times}$ and $y^{l_0}$.  

\begin{thm}[{{\cite[Th\'eoreme 5.2]{S3}}}]
For a maximal simple type $(J, \lambda)$, there exists an extension $\tilde{\lambda}$ of $\lambda$ to $\tilde{J}$.  
If $\tilde{\lambda}$ is such as an extension, $\cInd_{\tilde{J}}^{G} \tilde{\lambda}$ is irreducible and supercuspidal.  
\end{thm}

\subsection{Some property for $\Afra(E)$}

In this subsection, we fix an arbitrary finite field extension $E/F$ and will consider some property for $\Afra(E)$.  
We put $d=(\dim_{F}D)^{1/2}$.  

\begin{prop}
\label{periodofAfraE}
Then we have $\varpi_{F}\Afra(E)=\Pfra(E)^{de'/(d,e')}$, where $e'=e(E/F)$ is the ramification index of $E/F$.  
\end{prop}

\begin{prf}
Let $W$ be a simple right $E \otimes_{F} D$-module and $\Lambda = ( \Lambda_i )$ be an $\ofra_D$-chain in $W$ such that $\Afra(E)=\Afra(\Lambda)$.  
Applying \cite[Th\'eoreme 1.7]{SS} with $\Gamma = (\pfra_{D'}^i)$, we obtain $\rho_0 \in \N$ such that $\Pfra(E)^j \cap D' = \pfra_{D'}^{\lceil j/\rho_0 \rceil}$ for any $j \in \N$, where $\lceil \cdot \rceil$ is the ceil function.  
We have $\rho_0=d/(d, e'r_0)$ by the proof in \cite[Th\'eoreme 1.7]{SS}, where $r_0$ is the integer satisfying $\pfra_{D'}^i \varpi_{E}=\pfra_{D'}^{i+r_0}$, that is, $r_0=d/(d,[E:F])$.  
Then we obtain $\rho_0=(d, [E:F])/(d,e')$.  
Therefore, $v_{\Afra(E)}(\varpi_{D'})=\rho_0=(d,[E:F])/(d,e')$ and
\[
	v_{\Afra(E)}(\varpi_{F}) = e' v_{\Afra(E)}(\varpi_{E}) = e' \cdot \frac{d}{(d,[E:F])} v_{\Afra(E)}(\varpi_{D'}) = \frac{de'}{(d,e')}.  
\]
This implies that $\varpi_{F}\Afra(E)=\Pfra(E)^{de'/(d,e')}$.  
\end{prf}

\begin{cor}
\label{unramax}
Let $[\Afra, n, 0, \beta]$ is a simple stratum for some maximal simple type such that $E=F[\beta]$ is an unramified extension of $F$.  
Then $\Afra$ is a maximal hereditary $\ofra_F$-order.  
\end{cor}

\begin{prf}
Let $\Lambda=(\Lambda_i)$ be as in the proof of Proposition \ref{periodofAfraE}.  
Since $E/F$ is unramified, $\Lambda_{i} \pfra_{D}^d = \Lambda_{i} \varpi_{F} = \Lambda_{i+d}$ by Proposition \ref{periodofAfraE}.  
Then $\Lambda_{i} \pfra_{D} = \Lambda_{i+1}$.  
Therefore $\Afra(E)=\Afra(\Lambda)$ is simply equal to $\End_{\ofra_D}(\Lambda_0)$ and $\Afra(E)$ is maximal.  
Since $[\Afra, n, 0, \beta]$ is a simple stratum for some maximal simple type, $\Bfra = \Afra \cap \Cent_{A}(E)$ is a maximal $\ofra_E$-order in $\Cent_{A}(E)$.  
Then $\Afra$ is proper by Remark \ref{maxproper} (ii), and we have $\Afra \cong \M_{s}(\Afra(E))$ for some $s$ by Remark \ref{maxproper} (i).  
Therefore $\Afra$ is also maximal.  
\end{prf}

For later use, we will show the following lemma.  

\begin{lem}
\label{fundhered}
Let $E/F$ be an unramified extension and let $E_1$ be the unique subextension of $E/F$ such that $\dim_{F}E_1=(d, [E:F])$.  
We fix an inclusion $E_1 \hookrightarrow D$.  
Let $W$ be a simple right $E \otimes_{F} D$-module.  
Then there exists a right $D$-basis $\mathcal{B}$ of $W$ such that under the identification of $A(E)$ with $M_{m_0}(D)$ given by $\mathcal{B}$, 
\begin{enumerate}
\item we have $\Afra(E) = \M_{m_0}(\ofra_{D})$, and
\item the canonical inclusion $E_1 \subset E \hookrightarrow A(E)$ coincides with the composition of the fixed inclusion $E_1 \hookrightarrow D$ and the diagonal embedding of $D$ in $\M_{m_0}(D)$, 
\end{enumerate}
where $m_0=\dim_{D}(W)$.  
\end{lem}

\begin{prf}
First, we give an isomorphism of $W$ with a more concrete module.  
There is a (non-canonical) $E_1$-algebra isomorphism $E_1 \otimes _{F} D \cong \M_{m_1} (D_1)$ for some division $E_1$-algebra $D_1$ with $m_1 = [E_1 \colon F]$.  
We put $\dim_{E_1}{D_1} = d_1^2$ , and then $ d_1 = \left( d/[E_1 \colon F] \right)$.  
Since $d_1$ and $[E \colon E_1]$ are coprime by the assumption, $D'=E \otimes _{E_1} D_1$ is a division $E$-algebra and 
\[
	E \otimes _{F} D \cong E \otimes _{E_1} \left( E_1 \otimes _{F} D \right) \cong E \otimes _{E_1} \M_{m_1} (D_1) \cong \M_{m_1} \left( E \otimes_{E_1} D_1 \right) = \M_{m_1} \left( D' \right).  
\]
Then, we have $\dim _{D'} W = m_1 = [E_1 \colon F]$.  

Put $W_0=E \otimes _{E_{1}} D$.  
We show $W \cong W_0$.  
The module $W_0$ has a natural right $E \otimes _{F} D$-module structure. 
Then, there exists a right $E \otimes_{F} D$-module isomorphism $W_{0} \cong W ^{\oplus i}$ for some positive integer $i$.  
Since, in particular, this isomorphism is also an $E$-vector space isomorphism, the equation $\dim_{E}W_{0} = i \dim_{E}W$ holds.  
We have $\dim_{E}W_0 = \dim_{E_1}D = (\dim_{F}D)/[E_1 \colon F] = d^2 / [ E_1 \colon F]$ and $\dim_{E} W = d_1^2 \dim _{D'} W = d^2/[E_1 \colon F]$.  
Therefore, we obtain $i = 1$ and $W_0 \cong W$.  
Then we may assume $W = W_0$.  

We have $m_0 = \dim_{D}W_0 = [E : E_1]$.  
If an $E_1$-basis of $E$ is fixed, the embedding $E= \End_{E}(E) \subset \End_{E_1}(E) \cong \M_{m_0}(E_1)$ is defined.  
We fix an $\ofra_{E_1}$-basis $\mathcal{B}'=\{ b_1, \ldots, b_{m_0} \}$ of $\ofra_E$, and obtain a $D$-basis $\mathcal{B}=\{ b_1 \otimes 1, \ldots , b_{m_0} \otimes 1 \}$ of $W_0$.  
Under the idenfitication of $\M_{m_0}(D)$ with $A(E)$ defined by this basis, $E_1$ is embedded diagonally and $\ofra_{E}^{\times} \subset \GL_{m_0}(\ofra_{E_1}) \subset \GL_{m_0} (\ofra_{D})$.  

Put $\Afra_{0} = \M_{m_0}(\ofra_{D})$.  
Then, $\ofra_{E}^{\times}$ and $E_1 ^{\times} \hookrightarrow D^{\times}$ are contained in $\Kfra(\Afra_{0})$.  
Since $E / E_1$ is unramified, this shows that $\Afra_{0}$ is $E$-pure.  
Since $\Afra(E)$ is the unique $E$-pure hereditary $\ofra_{D}$-order, we have $\Afra(E) = \Afra_{0}$ and we obtain the desired conditions.  
\end{prf}

\section{Mackey Decomposition}
\label{decomp}

In this section, we will show that for $K=\GL_{m}(\ofra_D)$ and an irreducible supercuspidal representation $\pi$ of $G$, there exists a $[G, \pi]_{G}$-type $(K, \tau)$.  
Any maximal compact subgroup $K'$ of $G$ is $G$-conjugate to $K$, so this fact implies that there also exists a representation $\tau'$ of $K'$ such that $(K', \tau')$ is also a $[G, \pi]_{G}$-type.  

\begin{lem}
\label{ind_of_type}
Let $\mathfrak{s} \in \mathcal{B}(G)$, $(J, \lambda)$ be an $\mathfrak{s}$-type and $K$ be a compact open subgroup of $G$ containing $J$.  
Assume $\mathfrak{s}$ contains an irreducible supercuspidal representation $\pi_{0}$ of $G$.  
Then, the irreducible components of $\Ind_{J}^{K}\lambda$ contained in $\pi_{0}|_{K}$ are $\mathfrak{s}$-types.  
In particular, if $\Ind_{J}^{K}\lambda$ is irreducible, $\Ind_{J}^{K}\lambda$ is an $\mathfrak{s}$-type.  
\end{lem}

\begin{prf}
Let $\tau$ be an irreducible component of $\Ind_{J}^{K}\lambda$ contained in $\pi_{0}|_{K}$ and $\pi \in \Irr(G)$.  
We show that $\tau \subset \pi|_{K}$ if and only if $\pi \in \mathfrak{s}$.  

Let $\tau \subset \pi|_{K}$.  
Then, we have
\[
\left< \Ind_{J}^{K}\lambda, \pi \right> _{K} \geq \left< \tau, \pi \right>_{K} > 0, 
\]
and by the Frobenius reciprocity we obtain $\left< \lambda, \pi \right>_{J} >0$.  
Since $(J, \lambda)$ is an $\mathfrak{s}$-type, it follows that $\pi \in \mathfrak{s}$.  

Conversely, let $\pi \in \mathfrak{s}$.  
Then, there exists an unramified character $\chi$ of $F^{\times}$ with $\pi \cong \pi_{0} \otimes \chi \circ \Nrd_{A/F}$ and we have $\pi |_{K} \cong \pi_{0} |_{K}$.  
Therefore $\tau$ is a $K$-subrepresentation of $\pi |_{K} \cong \pi_{0} |_{K}$.  
\end{prf}

Let $\pi$ be an irreducible supercuspidal representations of $G$.  
There exists a simple type $(J, \lambda)$ and a unique extension $(\tilde{J}, \tilde{\lambda} )$ such that 

\[
	\pi = \cInd _{\tilde{J}} ^{G} \tilde{\lambda}.  
\]

If $(J, \lambda)$ is of level 0, then $J$ is the multiplicative group of some maximal hereditary $\ofra_F$-order $\Afra$ in $A$.  
If $(J,\lambda)$ attached to the simple stratum $[ \Afra, n, 0, \beta]$ with $n>0$.  
Since $\Bfra$ is maximal, $\Afra$ is a proper $F[\beta]$-pure hereditary $\ofra_F$-order.  
Anyway, by $G$-conjugation, we may assume $\Afra$ is contained in $\M_{m}(\ofra_D)$ and $\U(\Afra)$ is a group consisting of block upper triangular matrices modulo $\pfra_{D}$.  
Let $\Lambda=(\Lambda_i)$ be an $\ofra_D$-chain which $\Afra$ is the hereditary $\ofra_F$-order associated with, and let $e$ be the period of $\Lambda$.  
Then there is a natural bijection
\[
	\Afra/\Pfra \cong \prod_{i=0}^{e-1} \End_{k_{D}}(\Lambda_{i}/\Lambda_{i+1}).  
\]
We can take an $\ofra_D$-chain $\Lambda$ such that $\Afra = \Afra(\Lambda)$ and a basis $\{ v_1, \ldots , v_m \}$ of $V$ as
\begin{eqnarray*}
	\Lambda_0 & = & v_1 \ofra_{D} + \cdots + v_m \ofra_{D}, \\
	\Lambda_i & = & v_1 \ofra_{D} + \cdots + v_{(m/e)(e-i)} \ofra_{D} + v_{(m/e)(e-i)+1} \pfra_{D} + \cdots + v_{m} \pfra_{D}
\end{eqnarray*}
for $0<i<e$.  

\begin{lem}
\label{tildeJinKA}
If $(\tilde{J}, \tilde{\lambda})$ is as above, then $\tilde{J}$ is contained in $\Kfra(\Afra)$.  
\end{lem}

\begin{prf}
The group $\tilde{J}$ is generated by $J, F^{\times}$ and $y^{l_0}$ as in \S \ref{max_simple_types}.  
Since $J \subset \Afra$ and $F^{\times}$ is the center of $G$, it is enough to show that $y$ normalizes $\Afra$.  
Moreover, if $(J, \lambda)$ is of level 0, then $y \in \Kfra(\Afra)$ by definition of $y$.  
Then we may assume $(J, \lambda)$ is of level $>0$.  

Because $\Afra$ is proper and $(J, \lambda)$ is maximal, we may assume $\Afra = \M_{m'}\left( \Afra(E) \right)$.  
We fix a uniformizer $\varpi_{D'}$ of $D'$, contained in $A(E)$, and $y$ is a diagonal matrix such that its diagonal coefficients are all $\varpi_{D'}$.  
Then $y\Afra y^{-1}=\M_{m'}(\varpi_{D'}\Afra(E)\varpi_{D'}^{-1})$.  
Therefore, it is enough to show $\varpi_{D'}\Afra(E)\varpi_{D'}^{-1}=\Afra(E)$.  

Since $\varpi_{D'} \in D' = \Cent_{A(E)}(E)$ and $\Afra(E)$ is $E$-pure, for every $x \in E^{\times}$, we have
\[
x\left( \varpi_{D'}\Afra(E)\varpi_{D'}^{-1} \right) x^{-1} = \varpi_{D'} x \Afra(E) x^{-1} \varpi_{D'}^{-1} = \varpi_{D'}\Afra(E)\varpi_{D'}^{-1}.  
\]
Therefore it follows that $\varpi_{D'}\Afra(E)\varpi_{D'}^{-1}$ is $E$-pure.  
Since $\Afra(E)$ is the unique $E$-pure hereditary $\ofra_F$-order in $A(E)$, we obtain $\varpi_{D'}\Afra(E)\varpi_{D'}^{-1}=\Afra(E)$.  
\end{prf}

\begin{lem}
\label{cptint}
If $K_{0}$ is a compact open subgroup of $G$ containing $J$, then $\tilde{J} \cap K_{0}=J$.  
In particular, $J$ is the unique maximal compact subgroup of $\tilde{J}$.  
\end{lem}

\begin{prf}
If $(J, \lambda)$ is of level 0, then $J$ is a maximal compact subgroup of $G$, whence $J \subset K_{0}$ implies $J=K_{0}$.  

If $(J, \lambda)$ is of level $>0$, let $L$ be a subgroup of $\Kfra(\Bfra)$, generated by $F^{\times}$ and $y^{l_0}$.  
Since $y$ normalizes $J$, we have $\tilde{J} = LJ$ and $\tilde{J} \cap K_{0} = (L \cap K_{0})J$.  
The group $L \cap K_{0}$ is compact, whence contained in $\U(\Bfra)$.  
Since $\U(\Bfra)$ is a subgroup of $J$, we have $\tilde{J} \cap K_{0} = J$.  
\end{prf}

Put $\tilde{\rho}=\cInd_{\tilde{J}}^{\Kfra(\Afra)}\tilde{\lambda}$ and $\rho = \Ind_{J}^{\U(\Afra)}\lambda$.  
By the transivity of compact induction, we have $\pi \cong \cInd _{\Kfra(\Afra)}^{G}\tilde{\rho}$.  
By Lemma \ref{cptint}, we have $I_{G}(\lambda, \lambda) \cap \U(\Afra) = J$ and $\rho$ is irreducible.  

By the Mackey decomposition, we have 
\[
\Res_{K}^{G} \pi \cong \bigoplus _{g \in K \backslash G / \Kfra(\Afra)} \Ind _{ K \cap ^{g}\Kfra(\Afra)} ^{K} \Res _ { K \cap ^{g}\Kfra(\Afra)} ^{ ^{g}\Kfra(\Afra)} {^g}\tilde{\rho}.  
\]

By the term corresponding to g=1 in the above decomposition, we have $\pi |_{K} \supset \Ind _{\U (\Afra)}^{K} \Res_{\U (\Afra)}^{\Kfra(\Afra)} \tilde{\rho}$ since $K \cap \Kfra(\Afra) = \U (\Afra)$.  
By the Mackey decomposition, we obtain 
\[
	 \Res_{\U (\Afra)}^{\Kfra(\Afra)} \tilde{\rho} \cong \bigoplus _{h \in \U(\Afra) \backslash \Kfra(\Afra) / \tilde{J} } \Ind _{\U(\Afra) \cap ^{h}\tilde{J}} ^{\U(\Afra)} \Res_{\U(\Afra) \cap ^{h}\tilde{J} } ^{ ^{h}\tilde{J}} { ^{h}\tilde{\lambda}}.  
\]
The group $\U(\Afra)$ is a normal subgroup of $\Kfra(\Afra)$, hence $\U(\Afra) = {^{h}\U(\Afra)}$, $\U(\Afra) \cap {^{h}\tilde{J}} = {^{h} \left( \U(\Afra)\cap\tilde{J} \right) } = {^{h}J}$ and $\Res_{\U(\Afra) \cap ^{h}\tilde{J}}^{ ^{h}\tilde{J}} { ^{h}\tilde{\lambda}} ={ ^{h}\lambda}$.  
Therefore we obtain
\[
\Res_{\U (\Afra)}^{\Kfra(\Afra)} \tilde{\rho} \cong \bigoplus _{h \in \U(\Afra) \backslash \Kfra(\Afra) / \tilde{J} } { ^{h} \left( \Ind_{J}^{\U(\Afra)} \lambda \right) } \cong \bigoplus _{h \in \U(\Afra) \backslash \Kfra(\Afra) / \tilde{J} } { ^{h} \rho}.  
\]

\begin{cor}
\label{g1type}
Let $h \in \Kfra(\Afra)$.  
Then, $\Ind _{\U(\Afra)} ^{K} { ^{h} \rho}$ is irreducible and a $[G, \pi]_{G}$-type.  
\end{cor}

\begin{prf}
We have $\Ind _{\U(\Afra)} ^{K} { ^{h} \rho} = \Ind _{\U(\Afra)} ^{K} \Ind_{{^h}J} ^{\U(\Afra)} {^h}\lambda = \Ind _{{^h}J} ^{K} {^h}\lambda$.  
Since $({^h}J, {^h}\lambda)$ is a simple type, it is enough to show that $\Ind _{{^h}J} ^{K} {^h}\lambda$ is irreducible by Lemma \ref{ind_of_type}.  
We have
\[
	I_{K}({^h}\lambda, {^h}\lambda)=I_{G}({^h}\lambda, {^h}\lambda) \cap K ={^h}\tilde{J} \cap K = \tilde{({^h}J)} \cap K = {^h}J, 
\]
from $K \supset \U(\Afra) = {^h}\U(\Afra) \supset {^h}J$ and Lemma \ref{cptint}.  
Therefore, $\Ind _{{^h}J} ^{K} \lambda$ is irreducible.  
\end{prf}

\begin{prop}
For every $h \in \Kfra(\Afra)$, there exists $k \in N_{G}(K)$ such that $\Ind_{J}^{K} \lambda \cong {^{k}}(\Ind_{^{h}J}^{K} {^h}\lambda)$.  
\end{prop}

\begin{prf}
Assume that we find $k \in N_{G}(K)$ such that $kh \in \tilde{J}$.  
Then $kh$ is the element of $I_{G}(\lambda, \lambda)$ and
\begin{eqnarray*}
	0 & \neq & \Hom_{J \cap {^{kh}}J} \left( \Res_{J \cap {^{kh}}J}^{J} \lambda, \Res_{J \cap {^{kh}}J}^{{^{kh}}J}{^{kh}\lambda} \right) \cong \Hom_{J} \left( \lambda, \Ind_{J \cap {^{kh}}J}^{J} \Res_{J \cap {^{kh}}J}^{{^{kh}}J}{^{kh}\lambda} \right) \\
	& \hookrightarrow & \Hom_{J} \left( \lambda, \bigoplus_{j \in J \backslash K /{^{kh}J}} \Ind_{J \cap {^{jkh}}J}^{J} \Res_{J \cap {^{jkh}}J}^{{^{jkh}}J}{^{kh}\lambda} \right) \\
	& \cong & \Hom_{J} \left( \lambda, \Res_{J}^{K} \Ind_{{^{kh}}J}^{K} {^{kh}} \lambda \right) \cong \Hom_{K} \left( \Ind_{J}^{K} \lambda, {^k} \left( \Ind_{{^h}J}^{K}{^h} \lambda \right) \right).  
\end{eqnarray*}
Since $\Ind_{J}^{K}\lambda$ and $\Ind_{{^h}J}^{K}{^h} \lambda$ are irreducible by Corollary \ref{g1type}, they are isomorphic.  
Therefore, it is enough to show that there exists $k \in N_{G}(K)$ such that $kh \in \tilde{J}$.  

First, we fix a uniformizer $\varpi_D$ of $D$ and embed $\ofra_D$ in $\M_{m}(\ofra_D)$ diagonally.  
Then $\varpi_D \in \Kfra(\Afra)$ and $N_{G}(K)$ is generated by $K$ and $\varpi_D$.  

We will show that $v_{\Afra}\left( N_{G}(K) \cap \Kfra(\Afra) \right)$ and $v_{\Afra}(\tilde{J})$ generate $\Z$.  
If $(J, \lambda)$ is of level 0, then $\Afra$ is maximal and $v_{\Afra}(\varpi_{D}) = 1$, whence $v_{\Afra}\left( (N_{G}(K) \cap \Kfra(\Afra) \right)=\Z$.  
Then we may assume that $(J, \lambda)$ is of level $>0$ with a simple stratum $[\Afra, n, 0, \beta]$.  
We put $E=F[\beta]$.  
Since $\Afra$ is proper and $\Bfra$ is maximal, the period of $\Afra$ is equal to the period of $\Afra(E)$ and we have $v_{\Afra}(\varpi_{D}) = (1/d)v_{\Afra}(\varpi_{F}) = e'/(d,e')$ by Proposition \ref{periodofAfraE}.  
Hence, we have
\[
	v_{\Afra}\left( N_{G}(K) \cap \Kfra(\Afra) \right) \supset v_{\Afra}(\varpi_{D})\Z = \frac{e'}{(d,e')}\Z.  
\]
On the other hand, we choose $y \in B^{\times}$ and $l_0 \in \N_{>0}$ as in \ref{max_simple_types}.  
Then we have $v_{\Afra}(y)=v_{\Afra(E)}(\varpi_{D'})=(d,[E:F])/(d,e')$ by the proof of Proposition \ref{periodofAfraE}.  
Moreover, $l_0$ is a divisor of $|\Gal(k_{D'}/k_{E})|=d/(d,[E:F])$ by definition of $l_0$.  
Therefore
\[
	v_{\Afra}(\tilde{J}) \supset v_{\Afra}(y^{l_0})\Z = \frac{(d,[E:F])l_0}{(d,e')}\Z \supset \frac{d}{(d,e')}\Z
\]
and we obtain
\[
	v_{\Afra}\left( N_{G}(K) \cap \Kfra(\Afra) \right) + v_{\Afra}(\tilde{J}) \supset \frac{e'}{(d,e')}\Z + \frac{d}{(d,e')}\Z = \frac{(d,e')}{(d,e')}\Z = \Z.  
\]

For fixed $h$, by the above discussion, there exist $k_1 \in N_{G}(K) \cap \Kfra(\Afra)$ and $j \in \tilde{J}$ such that $v_{\Afra}(k_1)+v_{\Afra}(j)=v_{\Afra}(h^{-1})$.  
Then we have $v_{\Afra}(k_1 h j)=0$ and $k_2 = k_1 h j \in \U(\Afra)$.  
We put $k={k_2}^{-1}k_1 \in N_{G}(K)$, whence we obtain $kh={k_2}^{-1}k_1 h = j^{-1} \in \tilde{J}$ and complete the proof.  
\end{prf}

For a general $g \in G$,  
\[
\Res _ { K \cap ^{g}\Kfra(\Afra)} ^{ ^{g}\Kfra(\Afra)}{^g}\tilde{\rho} = \Res _ { K \cap ^{g}\U(\Afra)} ^{ ^{g}\U(\Afra)} \Res _ {{^g}\U(\Afra)} ^{ ^{g}\Kfra(\Afra)}{^g}\tilde{\rho} = \bigoplus _{h \in \U(\Afra) \backslash \Kfra(\Afra) / \tilde{J} } \Res _ { K \cap ^{g}\U(\Afra)} ^{ ^{g}\U(\Afra)} {^{gh}}\rho.  
\]
Hence, we obtain
\[
	\Res_{K}^{G} \pi \cong \bigoplus _{g \in K \backslash G / \Kfra(\Afra)} \bigoplus _{h \in \U(\Afra) \backslash \Kfra(\Afra) / \tilde{J} } \Ind_{K \cap {^g}\U(\Afra)}^{K} \Res _ { K \cap ^{g}\U(\Afra)} ^{ ^{g}\U(\Afra)} {^{gh}}\rho.  
\]

\begin{lem}
\label{keyprop0}
Let $g \in G$ and $h \in \Kfra(\Afra)$.  
Let $\tau$ be a representation of $K$ such that
\[
\left< \tau, \Ind_{K \cap {^g}\U(\Afra)}^{K}{^{gh}}\rho \right>_{K} \neq 0.  
\]
Then there exists $h' \in \Kfra(\Afra)$ such that
\[
	^{h}\rho = \Ind_{J'}^{\U(\Afra)} \lambda', \text{ and } \left< \tau, \Ind_{K \cap {^g}J'}^{K}{^g}\lambda' \right>_{K} \neq 0, 
\]
where $(J', \lambda')=({^{h'}}J, {^{h'}}\lambda)$ is a simple type.  
\end{lem}

\begin{prf}
By the Mackey decomposition, 
\begin{eqnarray*}
	& & \Res_{{^{g^{-1}}}K \cap \U(\Afra)}^{\U(\Afra)} {^{h}}\rho \\
	& \cong & \bigoplus _{u \in {^{g^{-1}}}K \cap \U(\Afra) \backslash \U(\Afra) /{^h}J} \Ind_{{^{g^{-1}}}K \cap \U(\Afra) \cap {^{uh}}J}^{{^{g^{-1}}}K \cap \U(\Afra)} \Res _{{^{g^{-1}}}K \cap \U(\Afra) \cap {^{uh}}J} ^{{^{uh}}J} {^{uh}} \lambda \\
	& = & \bigoplus _{u \in {^{g^{-1}}}K \cap \U(\Afra) \backslash \U(\Afra) /{^h}J} \Ind_{{^{g^{-1}}}K \cap {^{uh}}J}^{{^{g^{-1}}}K \cap \U(\Afra)} \Res _{{^{g^{-1}}}K \cap {^{uh}}J} ^{{^{uh}}J} {^{uh}} \lambda.  \\
\end{eqnarray*}
Then
\begin{eqnarray*}
	& & \Ind_{{^{g^{-1}}}K \cap \U(\Afra)}^{^{g^{-1}}K} \Res_{{^{g^{-1}}}K \cap \U(\Afra)}^{\U(\Afra)} {^{h}}\rho \\
	& \cong & \bigoplus _{u \in {^{g^{-1}}}K \cap \U(\Afra) \backslash \U(\Afra) /{^h}J}  \Ind_{{^{g^{-1}}}K \cap {^{uh}}J}^{{^{g^{-1}}}K} \Res _{{^{g^{-1}}}K \cap {^{uh}}J} ^{{^{uh}}J} {^{uh}} \lambda.  
\end{eqnarray*}
Therefore there exists $u \in \U(\Afra)$ such that 
\[
\left< \tau,  \Ind_{K \cap {^{guh}}J}^{K} {^{guh}} \lambda \right>_{K} =  \left< {^{g^{-1}}} \tau, \Ind_{{^{g^{-1}}}K \cap {^{uh}}J}^{{^{g^{-1}}}K} {^{uh}} \lambda \right>_{{^{g^{-1}}}K} \neq 0.  
\]
Then $h'=uh$ and $(J', \lambda')=({^{uh}}J, {^{uh}}\lambda)$ satisfies the conditions of Lemma.  
\end{prf}

\begin{prop}
\label{keyprop}
Let $g \in G$ and $h \in \Kfra(\Afra)$.  
Let $\tau$ be an irreducible representation of $K$ such that
\[
\left< \tau, \Ind_{K \cap {^g}\U(\Afra)}^{K}{^{gh}}\rho \right>_{K} \neq 0. 
\] 
Let $(J, \lambda)$ be a simple type such that 
\[
	^{h}\rho = \Ind_{J}^{\U(\Afra)} \lambda \text{ and } \left< \tau, \Ind_{K \cap {^g}J}^{K}{^g}\lambda \right>_{K} \neq 0.  
\]
Suppose for every irreducible summand $\xi$ of $\lambda|_{J \cap {^{g^{-1}}K}}$, there exists an irreducible representation $\lambda'$ of $J$ such that 
\[
	\left< \xi, \lambda' \right>_{J \cap {^{g^{-1}}K}} \neq 0 \text{ and } I_{G}(\lambda, \lambda')=\emptyset.  
\]
Then, $\tau$ cannot be a type.  
\end{prop}

\begin{prf}
Assume that $\tau$ is a type.  
Let $\tilde{\tau}$ be an extension of $\tau$ to $F^{\times}K$.  
Then, by \cite[(5.2)]{BK2}, there are finitely many unramified characters $\chi_{i}$ of $F^{\times}$ such that
\[
	\cInd _{F^{\times}K}^{G} \tilde{\tau} = \bigoplus_{i} \pi \otimes \left( \chi_{i} \circ \Nrd_{A/F} \right).  
\]
Hence every irreducible subrepresentaion of
\[
	\Res_{J}^{G} \cInd _{F^{\times}K}^{G} \tilde{\tau} \cong \bigoplus_{g' \in J \backslash G / F^{\times}K} \Ind_{J \cap {^{g'}}F^{\times}K}^{J} \Res_{J \cap {^{g'}}F^{\times}K}^{{^{g'}}F^{\times}K} {^{g'}}\tilde{\tau}
\]
 is contained in $\pi|_{J}$.  
 In particular, $\pi|_{J}$ contains every irreducible subrepresentation of 
 \[
 \Ind_{J \cap {^{g^{-1}}}K}^{J} \Res_{J \cap {^{g^{-1}}}K}^{{^{g^{-1}}}K} {^{g^{-1}}}\tau.  
 \]
Since 
\[
\left< {^{g^{-1}}} \tau, \lambda \right>_{{^{g^{-1}}}K \cap J} = \left< \tau, {^{g}} \lambda \right> _{K \cap {^g}J} = \left< \tau, \Ind_{K \cap {^g}J}^{K}{^g}\lambda \right>_{K} \neq 0, 
\]
there exists an irreducible summand $\xi$ of $\lambda|_{J \cap {^{g^{-1}}K}}$ such that $\left< {^{g^{-1}}} \tau, \xi \right>_{J \cap {^{g^{-1}}}K} \neq 0$.  
By our assumption, there exists an irreducible representation $\lambda'$ of $J$ such that 
\[
	\left< \xi, \lambda' \right>_{J \cap {^{g^{-1}}K}} \neq 0, \text{ and } I_{G}(\lambda, \lambda')=\emptyset.  
\]
Then
\begin{eqnarray*}
	\left< \lambda', \Ind_{J \cap {^{g^{-1}}}K}^{J} \Res_{J \cap {^{g^{-1}}}K}^{{^{g^{-1}}}K} {^{g^{-1}}} \tau \right>_{J} & = & \left< \lambda', {^{g^{-1}}} \tau \right>_{J \cap {^{g^{-1}}}K} \\
	& \geq & \left< \xi, {^{g^{-1}}} \tau \right>_{J \cap {^{g^{-1}}}K} > 0, 
\end{eqnarray*}
whence $\lambda'$ is a $J$-subrepresentation of $\Ind_{J \cap {^{g^{-1}}}K}^{J} \Res_{J \cap {^{g^{-1}}}K}^{{^{g^{-1}}}K} {^{g^{-1}}} \tau$, so $\lambda'$ is contained in $\pi|_{J}$.  
On the other hand, there exists an extension $(\tilde{J}, \tilde{\lambda})$ of $(J, \lambda)$ such that $\pi \cong \cInd_{\tilde{J}}^{G} \tilde{\lambda}$.  
Then
\[
	\lambda' \subset \pi |_{J} \cong \bigoplus_{g' \in J \backslash G / \tilde{J}} \Ind_{J \cap {^{g'}}\tilde{J}}^{J} \Res_{J \cap {^{g'}}\tilde{J}}^{^{g'} \tilde{J}} {^{g'}}\tilde{\lambda}
\]
and there is $g' \in G$ such that 
\[
	\lambda' \subset \Ind_{J \cap { ^{g'}}\tilde{J}}^{J} \Res_{J \cap { ^{g'}\tilde{J}}}^{ ^{g'}\tilde{J}} { ^{g'}}\tilde{\lambda} = \Ind_{J \cap { ^{g'}J}}^{J} \Res_{J \cap { ^{g'}J}}^{ ^{g'}J} { ^{g'}} \lambda, 
\]
where since $J \cap \tilde{({ ^{g'}J})}$ is a compact subgroup of $\tilde{({^{g'}J})}$, the group $J \cap \tilde{({^{g'}J})}$ is a subgroup in ${^{g'}J}$ by Lemma \ref{cptint} and we have  $J \cap \tilde{({ ^{g'}J})} =  J \cap \tilde{({ ^{g'}J})} \cap {^{g'}J} = J \cap { ^{g'}J}$.  
Therefore $\left< \lambda', { ^{g'}}\lambda \right>_{J \cap { ^{g'}}J} \neq 0$ and $I_{G}(\lambda, \lambda') \neq \emptyset$, which is a contradiction.  
\end{prf}

\section{Representatives of $K \backslash G/\Kfra(\Afra)$}
\label{nice_repre}

Fix a uniformizer $\varpi_{D}$ of $D$.  
Put
\[
	R=\Set{
	\begin{pmatrix}
	\varpi_{D}^{a_1} & & \\
	 & \ddots & \\
	 & & \varpi_{D}^{a_m} \\
	\end{pmatrix}
	| a_1, \ldots , a_m \in \Z
	}.  
\]
For $a_1, \ldots , a_m \in \Z$, the diagonal matrix whose entries are $\varpi_D^{a_1}, \varpi_D^{a_2}, \ldots , \varpi_D^{a_m}$ is denoted by $a_{R}(a_1, \ldots, a_m)$.  

\begin{lem}
\label{pick_repre}
Let $g'$ be an element of $G$ such that $Kg'\Kfra(\Afra) \neq K\Kfra(\Afra)$. 
Then, there exists an element $g$ of $R$ such that its entries are $\varpi_D^{a_1}, \varpi_D^{a_2}, \ldots , \varpi_D^{a_m}$ ($a_1, \ldots , a_m \in \N$), $a_{(m/e)i+1} \geq \cdots \geq a_{(m/e)(i+1)}$ for $0 \leq i < e$, and one of the following conditions holds:
\begin{enumerate}
\item There exists $0 \leq i < e$ and $1 \leq j < m/e$ such that $a_{(m/e)i+j}>a_{(m/e)i+j+1}$.  
\item We have $a_{(m/e)i+1} = \cdots = a_{(m/e)(i+1)}$ for $0 \leq i < e$ and 
\begin{itemize}
\item $a_{1} \geq 2$, 
\item there exists $2 \leq j \leq e$ such that for $1 \leq i \leq e$ we have $a_{(m/e)i} > 0$ if and only if $i < j$.  
\end{itemize}
\end{enumerate}
\end{lem}

\begin{prf}
Let $\Afra_{0} \subset \Afra$ be a minimal hereditary $\ofra_F$-order which consists of upper triangular matrices modulo $\pfra_{D}$.  
The groups of permutation matrices, isomorphic to $\mathfrak{S}_{m}$, is naturally embedded in $G$.  
Let $\tilde{W}$ be the semidirect product of $\mathfrak{S}_{m}$ and $R$ in $G$.  
Then we have $G=\coprod_{w \in \tilde{W}} \U(\Afra_{0}) w \U(\Afra_{0})$ by \cite[\S 0.8]{GSZ}.  
Therefore we obtain $G=\bigcup_{w \in \tilde{W}} Kw \Kfra(\Afra) = \bigcup_{g \in R}Kg\Kfra(\Afra)$.  
Then we can pick $g''=a_{R}(a_1, \ldots, a_m) \in Kg'\Kfra(\Afra) \cap R$.  
Since $a_{R}(1, \ldots, 1) \in \Kfra(\Afra)$, we may assume $a_1, \ldots, a_m \in \N$ and there exists $1 \leq i \leq m$ such that $a_i = 0$.  
We have ${\mathfrak{S}_{m/e}}^{\times e} \subset \U(\Afra)$, so we may assume $a_{(m/e)i+1} \geq \cdots \geq a_{(m/e)(i+1)}$ for $0 \leq i < e$.  

If there exists $0 \leq i < e$ and $1 \leq j < m/e$ such that $a_{(m/e)i+j}>a_{(m/e)i+j+1}$, then $g=g''$ satisfies (i).  
Otherwise, we define
\[
	t = 
	\begin{pmatrix}
	 & \I_{m/e} & & \\
	 & & \ddots & \\
	 & & & I_{m/e} \\
	\I_{m/e} & & & \\
	\end{pmatrix}
\]
and
\[
	b_R(b_1, \ldots, b_{e}) = 
	\begin{pmatrix}
	\varpi_D^{b_1} \I_{m/e} & & \\
	 & \ddots & \\
	 & & \varpi_D^{b_{e}} \I_{m/e} \\
	\end{pmatrix}
	\in R, 
\]
where $b_1, \ldots, b_e \in \Z$.  
Then we may assume there exist $b_1, \ldots, b_e \in \N$ such that $g''=b_{R}(b_1, \ldots , b_e)$ and $\min \left\{ b_1, \ldots, b_m\right\}=0$.  
We put
\[
	\Pi = 
	\begin{pmatrix}
	 & \I_{m/e} & & \\
	 & & \ddots & \\
	 & & & \I_{m/e} \\
	\varpi_{D}\I_{m/e} & & & \\
	\end{pmatrix}
	\in \Kfra(\Afra), 
\]
where $\I_{m/e}$ is the $(m/e) \times (m/e)$ identity matrix.  
We define a bijective map $f:\tilde{W} \to \tilde{W}$ as 
\[
f(w)=t^{-1}w\Pi, \quad w \in \tilde{W}.  
\]
We have $f \left( b_{R}(b_1, \ldots , b_e) \right) =b_{R}(b_{e}+1, b_1, \ldots , b_{e-1})$ and $f(g'') \in Kg'\Kfra(\Afra)$.  

We put $i_1 = \min \Set{1 \leq i \leq e | b_i = 0}$ and $b_{R}(b_1', \ldots, b_e')=f^{e-i_1}(g'')$.  
If $b_1' \geq 2$, then $g=f^{e-i_1}(g'')$ satisfies (ii) with $j=e$.  
If $b_1' =1$, we put $i_2 = \max \Set{i|b_{i'}' = 1\text{ for }i' \leq i}$.  
If $i_2=e-1$, then $({f^{-1}})^{e-1} \circ f^{e-i_1}(g'') =\I_{m}$ and $Kg'\Kfra(\Afra) =K\Kfra(\Afra)$, which is a contradiction.  
Therefore we have $i_2 < e-1$.  
Then $g=\left(f^{-1}\right)^{i_2} \circ f^{e-i_1}(g'')$ satisfies (ii) with $j=e-i_2$.  
\end{prf}

\begin{defn}
Let $g, g' \in G$.  
We say that a pair $\left(g, Kg'\Kfra(\Afra) \right)$ has Property A if $Kg'\Kfra(\Afra) \neq K\Kfra(\Afra)$ and $g$ is a representative of $Kg'\Kfra(\Afra)$ as in Lemma \ref{pick_repre} and satisfying (i).  
\end{defn}

\begin{lem}
\label{property_for_max}
If $\Afra$ is maximal, then for every double coset $Kg'\Kfra(\Afra)$ not equal to $K\Kfra(\Afra)$, there is $g \in Kg'\Kfra(\Afra)$ such that the pair $\left( g, Kg'\Kfra(\Afra) \right)$ has Property A.  
\end{lem}

\begin{prf}
Since $\Afra$ is maximal, we have $e=1$.  
We pick $g''=a_{R}(a_1, \ldots , a_m) \in Kg'\Kfra(\Afra) \cap R$ as in the proof of Lemma \ref{pick_repre}.  
The equation $e=1$ implies $a_1 \geq \cdots \geq a_m \geq 0$.  
If $a_1 = \cdots = a_m$, then $g'' \in \Kfra(\Afra) \subset K\Kfra(\Afra)$, which is a contradiction.  
Therefore there exists $1 \leq j < m$ such that $a_j > a_{j+1}$, so $\left(g'', Kg'\Kfra(\Afra) \right)$ has Property A.  
\end{prf}

\begin{prop}
\label{proppropA}
Assume that a pair $\left(g, Kg'\Kfra(\Afra) \right)$ has Property A.  
Then there is $0 \leq i_0 < e$ such that the image of $\U(\Afra) \cap {^{g^{-1}}}K$ under the map
\[
	\U(\Afra) \to \U(\Afra)/\U^{1}(\Afra) \cong \prod_{i=0}^{e-1} \Aut_{k_{D}}(\Lambda_{i}/\Lambda_{i+1}) \to \Aut_{k_{D}}(\Lambda_{i_0}/\Lambda_{i_0+1})
\]
is contained in some proper parabolic subgroup of $\Aut_{k_{D}}(\Lambda_{i_0}/\Lambda_{i_0+1})$.  
\end{prop}

\begin{prf}
For $0 < i \leq e$, the morphism $\U(\Afra) \to \Aut_{k_D}(\Lambda_{e-i}/\Lambda_{e-i+1}) \cong \GL_{m/e}(k_D)$ is as follows:
\[
	\begin{pmatrix}
	A_{1,1} & \cdots & A_{1,e} \\
	\vdots & \ddots & \vdots \\
	A_{e,1} & \cdots & A_{e,e} \\ 
	\end{pmatrix}
	\mapsto A_{i,i} \text{ mod } \pfra_{D},  
\]
where $A_{j,j'} \in \M_{m/e}(\ofra_{D})$ for $1 \leq j, j' \leq e$.  
Since $g$ is as in Lemma \ref{pick_repre} and satisfying (i), there exist $0 \leq i_1 < e$ and $1 \leq j_1 < m/e$ such that $a_{(m/e)i_1+j_1} > a_{(m/e)i_1+j_1+1}$.  
If $h=(h_{ij}) \in \U(\Afra) \cap {^{g^{-1}}}K$, then $h_{ij} \in \pfra_{D}$ for $(m/e)i_1+j_1 < i \leq (m/e)(i_1+1)$ and $(m/e)i_1+1 \leq j \leq (m/e)i_1 + j_1$.  
We put $i_0 = e - (i_1+1)$.  
Then the image of $h=(h_{ij}) \in \U(\Afra) \cap {^{g^{-1}}}K$ under the map $\U(\Afra) \to \Aut_{k_D}(\Lambda_{i_0}/\Lambda_{i_0+1})$ is contained in a proper parabolic subgroup of $\Aut_{k_D}(\Lambda_{i_0}/\Lambda_{i_0+1})$.  
\end{prf}

\section{Double Cosets with Property A}
\label{propertyA}
\begin{defn}
Suppose $\left(g, Kg\Kfra(\Afra) \right)$ has Property A.  
We define the subgroup $\mathcal{K}(g)$ of $\U(\Afra)$ by
\[
	\mathcal{K}(g) = \left( \U(\Afra) \cap { ^{g^{-1}}K} \right) \U^{1}(\Afra).  
\]
\end{defn}

\begin{defn}
\label{def_of_suff_sm}{{\cite[Definition 6.2]{Pas}}}
Let $H$ be a subgroup of $\GL_{N}(\F_q)$ and put
\[
	S= \Set{ h \in \GL_{N}(\F _q) | 
	\begin{array}{l}
	\text{$\chi_{h}(X) = f(X)^{l}$ for some irreducible} \\
	\text{polynomial $f \in \F_q[X]$ and $l \in \N$}
	\end{array}}, 
\]
where $\chi_{h}$ is the characteristic polynomial of $h$.  

The group $H$ is sufficiently small if there exists a proper subfield $\F$ of $\F_{q^{N}}$ such that for every $h \in H \cap S$ the roots of $\chi_{h}$ are the elements of $\F$.  
\end{defn}

The lemma in the following are proved in \cite[Lemma 6.5]{Pas}.  

\begin{lem}
\label{suffpp}
Let $\F _q/\F _{q'}$ be an extension of finite fields and $V_0$ be a finite dimensional $\F _q$-vector space.  
Let
\[
	\iota \colon \End_{\F_q}(V_0) \to \End_{\F _{q'}}(V_0)
\]
be the natural embedding by scalar restriction.  
If $H$ is a subgroup of $\Aut_{\F _q}(V_0)$ such that $\iota(H)$ is contained in a proper parabolic subgroup of $\Aut_{\F _{q'}}(V_0)$, then $H$ is sufficiently small.  
\end{lem}

\begin{prop}
\label{exam_suff_sm}
Suppose that a pair $\left(g, Kg\Kfra(\Afra)\right)$ has Property A and that $E/F$ is unramified.  
The image of $\mathcal{K}(g) \cap \U(\Bfra)$ in $\U(\Bfra)/\U^1(\Bfra)$ is sufficiently small. 
\end{prop}

\begin{prf}
Let $E_1/F$ be the unique field extension in $E$ as in Lemma \ref{fundhered}.  
Put $B_1 = \Cent_{A}(E_1)$ and $\Bfra_1 = \Afra \cap B_1$.  
We also fix an inclusion $E_1 \in D$.  

We will show that $\Bfra_1$ is a maximal hereditary $\ofra_{E_1}$-order.  
Let $W$ be a simple right $E \otimes_{F} D$-module.  
By Lemma \ref{fundhered}, there is a right $D$-basis $\mathcal{B}$ of $W$ such that $\Afra(E)$ is identified with $\M_{[E:E_1]}(\ofra_D)$ under the identification $A(E) \cong \M_{[E:E_1]}(D)$ induced by $\mathcal{B}$.  
On the other hand, since $\Afra$ is proper, there is an $E \otimes_{F} D$-isomorphism $V \cong W^{\oplus m'}$ which induces the identification between $A$ and $\M_{m'}(A(E))$ such that $\Afra = \M_{m'}(\Afra(E))$.  
We define a $D$-basis $\mathcal{B}_1$ of $V$ by $\mathcal{B}$ and the $D$-isomorphism $V \cong W^{\oplus m'}$, and obtain the identification $A \cong \M_{m}(D)$.  
Under this identification, the subrings $\Afra, E$ and $E_1$ in $A$ correspond the subrings in $\M_{m}(D)$ as follows:  
\begin{itemize}
\item $\Afra  = \M_{m'}(\M_{[E:E_1]}(\ofra_{D}))$.  
\item $E \hookrightarrow A(E) \hookrightarrow \M_{m'} \left(A(E)\right)$, where $E \hookrightarrow A(E)$ is the canonical embedding and $A(E) \hookrightarrow M_{m'}(A(E))$ is the diagonal embedding.  
\item $E_1 \hookrightarrow D \hookrightarrow \M_{m}(D)$, where $E_1 \hookrightarrow D$ is the inclusion we fixed above and $D \hookrightarrow M_{m}(D)$ is the diagonal embedding.  
\end{itemize}
We put $D_1=\Cent_{D}(E_1)$, and then $B_1=\Cent_{A}(E_1)=\M_{m}(D_1)$, so we have $\Afra \cap B_1 = \M_{m}(\ofra_{D_1})$ and $\Bfra_1$ is maximal.   

Since $E_1/F$ is unramified, we have $k_{D_1} = k_{D}$ and $\U(\Bfra_1)/\U^{1}(\Bfra_1) \cong \U(\Afra)/\U^{1}(\Afra)$.  
Since $\mathcal{K}(g) \supset \U^1(\Afra)$, under the identification $\U(\Bfra_1)/\U^{1}(\Bfra_1) \cong \U(\Afra)/\U^{1}(\Afra)$, the image of $\U(\Bfra_1) \cap \mathcal{K}(g)$ in $\U(\Bfra_1)/\U^1(\Bfra_1)$ is identified with $\mathcal{K}(g)/\U^1(\Afra)$, which is equal to the image of $\U(\Afra) \cap {^{g^{-1}}}K$ in $\U(\Afra)/\U^1(\Afra)$ and contained in some proper parabolic subgroup of $\U(\Afra)/\U^1(\Afra)$ by Proposition \ref{proppropA}, as $\Afra$ is maximal.  
Therefore the image of $\U(\Bfra_1) \cap \mathcal{K}(g)$ in $\U(\Bfra_1)/\U^1(\Bfra_1)$ is contained in some proper parabolic subgroup of $\U(\Bfra_1)/\U^1(\Bfra_1)$.  
We also have $B = \Cent_{B_1}(E)$, $\Bfra_1$ is $E$-pure, and that $\Bfra = \Bfra_1 \cap B$.  
Hence, we may assume that $E_1 = F$.  

The hereditary $\ofra_{F}$-order $\Afra$ is maximal, and we have $\Lambda_{i+1} = \Lambda_{i} \varpi_{D}$.  
Since $\Afra$ is $E$-pure, $\Lambda$ is also an $\ofra_{E} \otimes _{\ofra_{F}} \ofra_D$-chain, where we have $\ofra_{E} \otimes _{\ofra_{F}} \ofra_D \cong \ofra_{D'}$ because $[E \colon F]$ and $d$ are coprime and $E/F$ is unramified.  
Let $\Bfra'$ be the hereditary $\ofra_E$-order in $B=\End _{D'}(V)$ corresponding to $\Lambda$.  
Since $E/F$ is unramified, $\varpi_{D}$ is also a uniformizer of $D'$ and $\Lambda_{i+1} = \Lambda_{i}\pfra_{D'}$ holds.  
Therefore we have
\[
\U (\Bfra')/\U^{1} (\Bfra') \cong \Aut _{k_{D'}}(\Lambda_{0}/\Lambda_{1}).  
\]
On the other hand, $\Bfra' = \Set{ x \in B | x\Lambda_{i} \subset \Lambda_{i} , i \in \Z } = \Afra \cap B = \Bfra$ also holds.  
By Lemma \ref{suffpp} with $V_0 = \Lambda_{0}/\Lambda_{1}$, we obtain this proposition.  
\end{prf}

Next, we examine the restriction of some cuspidal representations of $\GL_{N}(\F_q)$ to its sufficiently small subgroup.  
Let $p$ be the characteristic of $\F_q$.  

\begin{lem}
\label{Zsigmondy}
For integers $q>1$ and $N>1$, there exists a prime number $r$ such that $r$ divides $q^N-1$, but does not divide $q^s-1$ for any $0<s<N$, except when $N=2$ and $q=2^i-1$ with $i \geq 2$ , or $N=6$ and $p=2$.  
\end{lem}

\begin{prf}
This result is known as Zsigmondy's Theorem in \cite{Zsig}.  
\end{prf}

\begin{prop}
\label{cusprepresonsuffsm}
Let $\sigma$ be an irreducible cuspidal representation of $\GL_{N}(\F _q)$ with character $\mathcal{X}$, and $H$ be a sufficiently small subgroup of $\GL_{N}(\F _q)$.  Suppose $[\F_q:\F_p] > 1$, and if $p=2$ and $[\F_q:\F_p] \leq 3$, we further assume that $H$ is contained in some proper parabolic subgroup of $\GL_{N}(\F _q)$.  
If $\xi$ is an irreducible representation of $H$ such that $\left<\xi, \sigma \right>_{H} \neq 0$, then there exists an irreducible representation $\sigma'$ of $\GL_{N}(\F _q)$ such that $\left< \xi, \sigma' \right>_{H} \neq 0$ and $\sigma' \neq \gamma \cdot \sigma$ for any $\gamma \in \Gal(\F _q/\F _p)$.  
\end{prop}

\begin{prf}
Let $S$ be as in Definition \ref{def_of_suff_sm}.  
Since $H$ is sufficiently small, there exists a proper subfield $\F$ of $\F_{q^N}$ such that for every $h \in H \cap S$ the roots of $\chi_h$ are in $\F$.  
Let $\Gamma = GL_{N}(\F_q), a=[\F _q : \F _p]$ and $b=[\F : \F_p]$.  

According to \cite{Gel}, for a character $\Psi$ of $\F_{q^{N}}^{\times}$ such that $\Psi^{q^s -1} \neq 1$ for any $s$ dividing but not equal to $N$, we can define a class function $\mathcal{X}_\Psi$ of $\Gamma$, which is a character of an irreducible cuspidal representation of $\Gamma$.  
Every character of an irreducible cuspidal representation of $\Gamma$ is obtained from a character of $\F_{q^{N}}^{\times}$ as above.  
If $\Theta$ is another character of $\F _{q^{N}}^{\times}$ such that $\Theta ^{q^s -1} \neq 1$ for any $s$ dividing but not equal to $N$, we have $\mathcal{X}_\Psi = \mathcal{X}_{\Theta}$ if and only if $\Theta = \Psi^{q^s}$ for some $s \geq 0$.  
Moreover, if $\Psi|_{\F^{\times}}=\Theta|_{\F^{\times}}$, then $\mathcal{X}_{\Psi}|_{H} = \mathcal{X}_{\Theta}|_{H}$.  

For every $\gamma \in \Gal(\F_q/\F_p)$, we take a character $\Psi_{\gamma}$ of $\F_{q^N}^{\times}$ corresponding to $\gamma \cdot \sigma$.  
Let $\Psi = \Psi_{1}$.  
Suppose there exists another character $\Theta$ of $\F_{q^{N}}^{\times}$ such that $\Theta ^{q^s -1} \neq 1$ for any $s$ dividing but not equal to $N$, $\Psi|_{\F^{\times}} = \Theta|_{\F^{\times}}$, and that $\Theta \neq \Psi_{\gamma}^{q^s}$ for any $s \geq 0$ and $\gamma \in \Gal(\F_q/\F_p)$.  
Let $\sigma'$ be a cuspidal representation of $\Gamma$ with character $\mathcal{X}_\Theta$.  
Then $\sigma' \neq \gamma \cdot \sigma$ for any $\gamma \in \Gal(\F_q/\F_p)$ because $\Theta \neq \Psi_{\gamma}^{q^s}$ for any $s$.  
Moreover, we have $\mathcal{X}_{\Psi}|_{H} = \mathcal{X}_{\Theta}|_{H}$, so $\sigma|_{H} = \sigma'|_{H}$.  

In almost all cases, we can pick such $\Theta$ by counting characters of $\F_{q^N}^{\times}$.  
By Lemma \ref{Zsigmondy}, there exists a prime number $r$ such that $r$ divides $p^{aN}-1=q^N-1$, but does not divide $p^s-1$ for any $0<s<aN$, except when $aN=2$ and $p=2^i-1$ with $i \geq 2$, or $aN=6$ and $p=2$.  

Suppose such $r$ as above exists.  
If $\Theta$ is not an $r$th power, $r$ divides the order of $\Theta$ since $\F_{q^N}^{\times}$ is cyclic.  
Then $\Theta^{q^s -1} \neq 1$ for any $0<s<N$ as $0<as<aN$ and $p^{as} -1 =q^{s} -1$ is not divided by $r$.  
In particular, $\Theta^{q^s -1} \neq 1$ for any $s$ dividing but not equal to $N$.  

We count the characters $\Theta$ such that $\Theta$ is not an $r$th power and $\Psi|_{\F^{\times}}=\Theta|_{\F^{\times}}$.  
At first, there exist $(q^{N}-1)/(p^b-1)$ characters such that $\Psi|_{\F^{\times}}=\Theta|_{\F^{\times}}$.  
For any character $\chi$ of $\F^{\times}$, let $\mathscr{C}(\chi)$ be the set of characters $\Theta$ of $\F_{q^N}^{\times}$ such that $\Theta$ is an $r$th power and $\Theta|_{\F^{\times}}=\chi$.  
Since $r$ does not divide $p^b-1 = |\F^{\times}|$, we can take a character $\Theta_0 \in \mathscr{C}(\Psi|_{\F^{\times}})$.  
By multiplying $\Theta_{0}^{-1}$, we obtain one-to-one correspondence between $\mathscr{C}(\Psi|_{\F^{\times}})$ and $\mathscr{C}(\mathbf{1})$.  
The characters of $\F_{q^N}^{\times}$ trivial on $\F^{\times}$ correspond to the characters of $\F_{q^N}^{\times} / \F^{\times}$.  
Under this correspondence, the elements of $\mathscr{C}(\mathbf{1})$ correspond to the characters that are $r$th power, since $r$ does not divide $|\F^{\times}|$.  
Therefore, there will be
\[
	\left( 1-\frac{1}{r} \right) \frac{q^{N}-1}{p^b-1}
\]
characters $\Theta$ such that $\Theta$ is not an $r$th power and $\Psi|_{\F^{\times}}=\Theta|_{\F^{\times}}$.  
The cardinality of the set $\left\{ \Psi_{\gamma}^{q^{s}} \mid \gamma \in \Gal(\F_q/\F_p) , s \in \N \right\}$ is at most $aN$, so if the inequation
\[
	\left( 1-\frac{1}{r} \right) \frac{q^{N}-1}{p^b-1} > aN
\]
holds, then we can take the desired character $\Theta$.  
Since $q^N=p^{aN}$ and $p^b = |\F| \leq |\F_{q^N}|^{1/2} = p^{aN/2}$, it is enough to show
\[
	\left( 1-\frac{1}{r} \right) (p^{aN/2}+1) > aN.  
\]
By our assumption, we have $a \geq 2$ and $aN \geq 2$.  
By using induction on $aN$, we can examine when the inequation holds.  
When $p \geq 5$, for any prime $r$ the inequation holds.  
When $p = 3$, we may assume $r \geq 5$ and then the inequation holds.  
When $p = 2$, since $3=2^2-1$, we may assume $r \geq 5$ and then the inequation holds for $aN = 3$ or $aN > 5$.  

The remaining cases are following:  $(aN, p)=(6,2),(4,2),(2,2)$, or $(2,2^i-1)$ with $i \geq 2$.  

When $aN=2$, we have $a=2$ and $N=1$, so $q=p^2$ and $\F = \F_p$.  
Since $N=1$, every irreducible representation of $\Gamma$ is a character.  
Our assumption that $H$ is sufficiently small implies $H \subset \F_p ^{\times}$. 
If $\chi$ is a character of $\F_p ^{\times}$, there exist $p+1$ characters of $\F_q ^{\times}$ extending $\chi$.  
Since $p+1 > 2$, there exists a character $\sigma'$ of $\F_q ^{\times}$ extending $\sigma|_{H}$ such that $\sigma' \neq \sigma, \gamma \cdot \sigma$, where $\gamma \in \Gal(\F_q/\F_p)$ is the nontrivial element.  

Suppose $(aN,p)=(6,2)$ or $(4,2)$.  

If $N=1$, then $b \leq a/2$ and $H \subset \F_{2^{b}}^{\times}$.  
If $\chi$ is a character of $\F_{2^{b}}^{\times}$, there exist at least $2^{a/2}+1$ characters of $\F_q ^{\times}$ extending $\chi$.  
Because $2^2+1>4$ and $2^3+1>6$, we can deduce the proposition for this case in the same way as in the case that $aN=2$.  

If $N>1$, then $a \leq 3$ and $H$ is contained in some proper parabolic subgroup $P$ of $\Gamma$ by our assumption.  
Let $U$ be the unipotent radical of the parabolic subgroup opposite to $P$.  
Let $\xi$ be as in this proposition.  
Since $U \cap P = 1$, we have $U \cap H = 1$.  
Hence, $\Res_{U}^{\Gamma} \Ind_{H}^{\Gamma} \xi \supset \Ind_{1}^{U} \Res_{1}^{H} \xi \supset \mathbf{1}_{U}$.  
Therefore there exists some irreducible summand $\sigma'$ of $\Ind_{H}^{\Gamma} \xi$ such that $\left< \mathbf{1}_{U} , \sigma' \right>_{U} \neq 0$.  
The inclusion $\sigma' \subset \Ind_{H}^{\Gamma} \xi$ implies that $\left< \xi, \sigma' \right>_{H} \neq 0$.  
Moreover, $\sigma'$ is not cuspidal, whence $\sigma' \neq \gamma \cdot \sigma$ as every $\gamma \cdot \sigma$ is cuspidal.  
\end{prf}

\begin{prop}
\label{lv0unicity}
Suppose $\Afra$ is a maximal hereditary $\ofra_{F}$-order of $A$.  
Let $\pi$ be an irreducible supercuspidal representation of $G$ with the $\U^1(\Afra)$-fixed part $\pi^{\U^1(\Afra)} \neq 0$.  
Then, there exists an irreducible cuspidal representation $\tau$ of $\U(\Afra)/\U^1(\Afra)$ such that $\pi |_{\U(\Afra)} $ contains the lift $\sigma$ of $\tau$ to $\U(\Afra)$ and for every irreducible $\U(\Afra)$-subrepresentation $\sigma'$ of $\pi^{\U^1(\Afra)} $, there exists $\gamma \in \Gal(k_D/ k_F)$ with $\sigma' \cong \gamma \cdot \sigma$.  
\end{prop}

\begin{prf}
This is the result from \cite[Theorem 5.5 (ii)]{GSZ}.  
\end{prf}

\begin{cor}
\label{cor1_of_lv0unicity}
Let $\Afra$ be a maximal hereditary $\ofra_F$-order of $A$ and $\sigma, \sigma'$ be irreducible representations of $\U(\Afra)$, trivial on $\U^1(\Afra)$.  
Suppose that $\sigma$ is cuspidal.  
Then, $\sigma$ and $\sigma'$ intertwine in G if and only if there exists $\gamma \in \Gal(k_D/ k_F)$ such that $\sigma' \cong \gamma \cdot \sigma$.  
\end{cor}

\begin{prf}
There exists $h \in \Kfra(\Afra)$ such that ${^h}\sigma \cong \gamma \cdot  \sigma$.  
Then, if $\sigma' \cong \gamma \cdot \sigma$, these representations $\sigma$ and $\sigma' \cong {^h}\sigma$ intertwine.  

Conversely, suppose that $\sigma$ and $\sigma'$ intertwine.  
Then, there exists $h \in G$ such that
\[
\Hom_{\U(\Afra) \cap {^{h}\U(\Afra)}} \left( \sigma', {^h}\sigma \right) \neq 0.  
\]
Let $(\tilde{J}, \tilde{\sigma})$ be an extension of the simple type $(\U(\Afra), \sigma)$ of level 0 such that $\pi = \cInd _{\tilde{J}} ^{G} \tilde{\sigma}$ is irreducible and supercuspidal.  
Since $\U(\Afra) \cap {^h}\tilde{J}$ is a compact subgroup of ${^h}\tilde{J}$, we have $\U(\Afra) \cap {^h}\tilde{J} \subset {^h}\U(\Afra)$ by Lemma \ref{cptint}.  
Then we obtain $\U(\Afra) \cap {^h}\tilde{J} = \U(\Afra) \cap {^h}\U(\Afra)$ and $^{h}\tilde{\sigma}|_{\U(\Afra) \cap {^h}\tilde{J}}={^{h}}\sigma|_{\U(\Afra) \cap {^h}\U(\Afra)}$, so we have
\begin{eqnarray*}
0 \neq \left< \sigma', {^h}\sigma \right>_{\U(\Afra) \cap {^{h}\U(\Afra)}} & = & \left<{^h}\tilde{\sigma}, \Res_{\U(\Afra) \cap {^{h}\tilde{J}}}^{\U(\Afra)} \sigma' \right>_{\U(\Afra) \cap {^{h}\tilde{J}}} \\
& = & \left< \Ind_{\U(\Afra) \cap {^{h}\tilde{J}}}^{\U(\Afra)} \Res_{\U(\Afra) \cap {^{h}\tilde{J}}}^{{^h}\tilde{J}}{^h}\tilde{\sigma} , \sigma' \right>_{\U(\Afra)} \\
& = & \dim \Hom _{\U(\Afra)}\left( \sigma', \Ind_{\U(\Afra) \cap {^{h}\tilde{J}}}^{\U(\Afra)} \Res_{\U(\Afra) \cap {^{h}\tilde{J}}}^{{^h}\tilde{J}}{^h}\tilde{\sigma} \right).  
\end{eqnarray*}
Since $\sigma'$ is irreducible, there exists an inclusion of $\U(\Afra)$-representations 
\[
	\sigma' \hookrightarrow \Ind_{\U(\Afra) \cap {^{h}\tilde{J}}}^{\U(\Afra)} \Res_{\U(\Afra) \cap {^{h}\tilde{J}}}^{{^h}\tilde{J}}{^h}\tilde{\sigma} \subset \Res_{\U(\Afra)}^{G} \cInd_{\tilde{J}}^{G} \tilde{\sigma} = \pi |_{\U(\Afra)}.  
\]
By Proposition \ref{lv0unicity}, there exists $\gamma \in \Gal(k_D/k_F)$ such that $\sigma' \cong \gamma \cdot \sigma$.  
\end{prf}

\begin{cor}
\label{key_to_PropertyA}
Let $\Afra$ be a maximal hereditary $\ofra_F$-order of $A$, and let $\mathcal{K}$ be a subgroup of $\U(\Afra)$ such that the image of $\mathcal{K}$ in $\U(\Afra)/\U^1(\Afra)$ is sufficiently small, and $\sigma$ be a lift of an irreducible cuspidal representation $\tau$ of $\U(\Afra)/\U^1(\Afra)$ to $\U(\Afra)$.  
Moreover, if $q_F=2$ and $q_D \leq 8$, we further assume that the image of $\mathcal{K}$ in $\U(\Afra)/\U^1(\Afra)$ is contained in some proper parabolic subgroup of $\U(\Afra)/\U^1(\Afra)$.  
Then, for every $\mathcal{K}$-subrepsentation $\zeta$ of $\sigma$, there exists an irreducible representation $\sigma'$ of $\U(\Afra)$ such that $\sigma'$ is trivial on $\U^1(\Afra)$, and we have $\left< \zeta , \sigma' \right>_{\mathcal{K}} \neq 0$ and $I_{G}\left( \sigma, \sigma' \right)=\emptyset$.  
\end{cor}

\begin{prf}
Let $H$ be the image of $\mathcal{K}$ in $\U(\Afra)/\U^{1}(\Afra)$.  
Let $\omega$ be a $H$-representation such that $\zeta$ is the lift of $\omega$.  

Suppose there is an irreducible representation $\tau'$ of $\U(\Afra)/\U^{1}(\Afra)$ such that $\left< \omega, \tau' \right> _{H} \neq 0$ and $\tau' \neq \gamma \cdot \tau$ for any $\gamma \in \Gal(k_{D}/k_{F})$.  
Let $\sigma'$ be the lift of $\tau'$ to $\U(\Afra)$.  
Then $\left< \zeta, \sigma' \right> _{\mathcal{K}} = \left< \omega, \tau' \right> _{H} \neq 0$ and $\sigma' \neq \gamma \cdot \sigma$ for any $\gamma$.  
By Corollary \ref{cor1_of_lv0unicity}, the latter condition implies $I_{G}\left( \sigma, \sigma'\right)=\emptyset$.  
Therefore, it is enough to show that there exists such $\tau'$.  

Unless the characteristic of $k_{F}$ is equal to 2, or $q_D \leq 8$, by Proposition \ref{cusprepresonsuffsm}, such a representation $\tau'$ exists.  
When $q_F=2$ and $q_D \leq 8$, since $H$ is contained in some proper parabolic subgroup, we can apply Proposition \ref{cusprepresonsuffsm} and obtain $\tau'$ with desired conditions.  
If $q_F>2$ and $q_D=4,8$, then $q_D=q_F$ and $D=F$.  
Since $q_F \geq 4$, this corollary for this case is already proved in \cite[Corollary 6.14]{Pas}.  
\end{prf}

\begin{prop}
\label{exam_intertwining}
Let $[\Afra, n, 0, \beta]$ be a simple stratum of $A$ with $n>0$, and $\theta \in \mathscr{C}(\beta, \Afra)$.  
Let $\kappa$ be a $\beta$-extension of the Heisenberg representation $\eta$ of $\theta$.  
Let $\sigma$ and $\sigma'$ be irreducible representations of $J$, trivial on $J^1$.  
If $I_{B^{\times}}(\sigma, \sigma')=\emptyset$, then $I_{G}(\kappa \otimes \sigma, \kappa \otimes \sigma') = \emptyset$.  
\end{prop}

\begin{prf}
Let $x$ be an element of $G$ and $\phi \in I_{x}(\kappa \otimes \sigma, \kappa \otimes \sigma')$.  
We will show that $\phi = 0$.  
Let $X, Y$ and $Y'$ be the representation spaces of $\kappa, \sigma$ and $\sigma'$.  
There exist $S_j \in \End_{\C}(X)$ and $T_j \in \Hom_{\C}(Y, Y')$ such that $\phi = \sum_{j}S_{j} \otimes T_{j}$ and $T_{j}$ are linearly independent.  
For $h \in J^{1} \cap {^x}J^1$, we have
\[
	\phi \circ \left( \eta(h) \otimes id_{Y} \right) = \phi \circ \left( ( \kappa \otimes \sigma )(h) \right) \\
	 = \left( {^x} ( \kappa \otimes \sigma' )(h) \right) \circ \phi \\
	 = \left( {^x} \eta (h) \otimes id_{Y'} \right) \circ \phi.  \\
\]
Since $\phi = \sum_{j}S_{j} \otimes T_{j}$, we also have
\begin{eqnarray*}
	\phi \circ \left( \eta(h) \otimes id_{Y} \right)  & = & \sum_{j} \left( S_{j} \circ \eta(h) \right) \otimes T_{j}, \\
	\left( {^x} \eta (h) \otimes id_{Y'} \right) \circ \phi & = & \sum_{j} \left( {^x} \eta (h) \circ S_{j} \right)  \otimes T_{j}, 
\end{eqnarray*}
and $\sum_{j} \left( S_{j} \circ \eta(h) - {^x} \eta (h) \circ S_{j} \right) \otimes T_{j} = 0$.  
Since $T_j$ are linearly independent, it follows that $S_{j} \circ \eta(h) = {^x} \eta (h) \circ S_{j}$, that is, $S_j \in I_{x}(\eta, \eta)$.  
Now $\dim _{\C} I_{x}(\eta, \eta) \leq 1$.  
Then we may assume $\phi = S \otimes T$ for some $S \in I_{x}(\eta, \eta) = I_{x}(\kappa, \kappa)$ and $T \in \Hom_{\C}(Y, Y')$.  
We may assume $S \neq 0$.  
Then $x \in I_{G}(\kappa, \kappa) = JB^{\times}J$, so we may assume $x \in B^{\times}$.  
For $h \in J \cap {^x}J$, we have
\[
	\left( {^x} \kappa (h) \circ S \right) \otimes \left( T \circ \sigma (h) \right) = \left( S \circ \kappa (h) \right)  \otimes \left( T \circ \sigma (h) \right) = \left( {^x} \kappa (h) \circ S \right) \otimes \left( {^x} \sigma ' (h) \otimes T \right).  
\]
Since ${^x} \kappa (h) \circ S \neq 0$, it follows that $T \circ \sigma (h) = {^x} \sigma ' (h) \circ T$, that is, $T \in I_{x}(\sigma, \sigma')$.  
The assumption $I_{B^{\times}}(\sigma, \sigma')=\emptyset$ implies $I_{x}(\sigma, \sigma') = 0$, whence $T=0$ and $\phi = 0$.  
\end{prf}

\begin{lem}
\label{lv0twist}
Let $[\Afra, n, 0, \beta]$ be a simple stratum such that $\beta \in F$ and $\Afra$ is a maximal $\ofra_{F}$-order in $A$.  
\begin{enumerate}
\item We have $J(\beta, \Afra)=\U(\Afra)$ and $J^1(\beta, \Afra)=H^1(\beta, \Afra)=\U^1(\Afra)$.  
\item Let $\theta \in \mathscr{C}(\beta, \Afra)$ and let $\kappa$ be a $\beta$-extension of $\theta$.  
Then $\kappa$ is the restriction of a character of $G$ to $U(\Afra)$.  
\item Let $(J, \lambda)$ be attached to $[\Afra, n, 0, \beta]$.  
Then $(J, \lambda)$ is obtained from twisting a simple type of $G$ of level 0 by some character of $G$.  
\end{enumerate}
\end{lem}

\begin{prf}
Since $\beta \in F$, we have $B=\Cent_{A}(\beta)=A$.  
Then $\Bfra=\Afra \cap B=\Afra$.  

To show (i), we recall some property for $J(\beta, \Afra)$ and $H(\beta, \Afra)$ from \cite{S1}.  
The subrings $\mathfrak{J}(\beta, \Afra)$ and $\mathfrak{H}(\beta, \Afra)$ of $A$ are defined in \cite[\S 3.3.1]{S1}.  
According to \cite[\S 3.3]{S1}, the group $J(\beta, \Afra)$ (resp. $H(\beta, \Afra)$ ) is the multiplicative group of $\mathfrak{J}(\beta, \Afra)$ (resp. $\mathfrak{H}(\beta, \Afra)$ ).  

Therefore, for (i), it is enough to show that $\mathfrak{J}(\beta, \Afra)=\mathfrak{H}(\beta, \Afra)=\Afra$.  
Since $\beta \in F$, the element $\beta$ is 'minimal' in the sense of \cite[2.3.3]{S1}.  
Then we have $\mathfrak{J}(\beta, \Afra)=\mathfrak{H}(\beta, \Afra)=\Afra$ by \cite[Proposition 3.42(1)]{S1}.  
  
The simple character $\theta$ factors through $\Nrd_{A/F}$ by \cite[Proposition 3.47(3)]{S1}.  
Then there exists a character $\chi_1$ of $F^{\times}$ such that $\theta = \chi_1 \circ \left( \Nrd_{A/F} | _{\U^1(\Afra)} \right)$.  
We put $\kappa_0=\chi_1 \circ \left( \Nrd_{A/F} | _{\U(\Afra)} \right)$.  
Then $\kappa_0$ is an extension of $\theta$ to $\U(\Afra)$.  

We will show that $\kappa_0$ is a $\beta$-extension of $\theta$.  
To show this, it is enough to show that $I_{G}(\kappa_0, \kappa_0)=G$.  
Since $\kappa_0$ is the restriction of a character of $G$, we have $\kappa_0(x^{-1}gx)=\kappa_0(g)$ for any $x \in G$ and $g \in \U(\Afra)$.  
Therefore $x \in I_{G}(\kappa_0, \kappa_0)$, that is, $I_{G}(\kappa_0, \kappa_0)=G$.  

There exists a character $\chi_2$ of $F^{\times}$ such that $\kappa =  \kappa_0 \otimes \left( (\chi_2 \circ \Nrd_{A/F})|_{\U(\Afra)} \right)$ by Th\'eoreme 2.28 in \cite{S2}.  
Since $\kappa_0$ and $(\chi_2 \circ \Nrd_{A/F}) |_{\U(\Afra)}$ are the restrictions of characters of $G$, $\kappa$ is also the restrition of a character of $G$.  

To show (iii), let $\kappa$ and $\sigma$ be irreducible $J(\beta, \Afra)$-representations such that $\lambda = \kappa \otimes \sigma$, as in Defenition \ref{def_of_simple_type}.  
Since $\Bfra=\Afra$, the representation $\sigma$ is a simple type of $G$ of level 0.  
Moreover, $\kappa$ is the restriction of a character of $G$ by (ii).  
Therefore, (iii) holds.  
\end{prf}

\begin{prop}
\label{conclude_PropertyA}
Suppose $\left( g, Kg'\Kfra(\Afra) \right)$ has Property A.  
Let $(J_0, \lambda_0)$ be a maximal simple type such that $J_0 \subset \U(\Afra)$.  
If $(J_0, \lambda_0)$ is of level $>0$, then we further assume that $(J_0, \lambda_0)$ is attached to a simple stratum $[\Afra, n, 0, \beta_0]$ such that $F[\beta_0]/F$ is unramified.  
We put $\rho = \Ind_{J_0}^{\U(\Afra)} \lambda_0$.  
Let $\tau$ be an irreducible $K$-representation such that
\[
\left< \tau, \Ind_{K \cap {^g}\U(\Afra))}^{K} {^{gh}} \rho \right>_{K} \neq 0
\]
for some $h \in \Kfra(\Afra)$.  
Then, $\tau$ cannot be a type.  
\end{prop}

\begin{prf}
By Lemma \ref{keyprop0}, there exists $h' \in \Kfra(\Afra)$ such that $(J, \lambda)=({^{h'}}J_0, {^{h'}}\lambda_{0})$ is a simple type, 
\[
	^{h}\rho = \Ind_{J}^{\U(\Afra)} \lambda \text{ and } \left< \tau, \Ind_{K \cap {^g}J}^{K}{^g}\lambda \right>_{K} \neq 0
\]
hold.  

Assume that $(J_{0}, \lambda_0)$ is level 0.  
Then $J_0=\U(\Afra)$ is maximal and $\lambda_0$ is trivial on $\U^1(\Afra)$.  
Since $h' \in \Kfra(\Afra)$, we have $J=\U(\Afra)$ and $\lambda$ is also trivial on $\U^1(\Afra)$.  
Therefore, the image of $\mathcal{K}(g)$ in $\U(\Afra)/\U^1(\Afra)$ is contained in some proper parabolic subgroup by applying Proposition \ref{proppropA} with $e=1$.  
Let $\xi$ be an irreducible $J \cap {^{g^{-1}}}K$-subrepresentation of $\lambda$.  
Since $\lambda$ is trivial on $\U^1(\Afra)$, $\xi$ can be extended to a unique $\mathcal{K}(g)$-subrepresentation of $\lambda$ by $\xi |_{\U^1(\Afra)}$ being trivial.  
Hence by Corollary \ref{key_to_PropertyA}, there exists an irreducible representation $\lambda'$ of $\U(\Afra)$ such that the representation $\lambda'$ is trivial on $\U^1(\Afra)$, $\left< \xi , \lambda' \right>_{\mathcal{K}(g)} \neq 0$ and $I_{G}\left( \lambda, \lambda' \right)=\emptyset$.  
Since $J \cap {^{g^{-1}}}K \subset \mathcal{K}(g)$, we have $\left< \xi , \lambda' \right>_{J \cap {^{g^{-1}}}K} \geq \left< \xi, \lambda' \right>_{\mathcal{K}(g)} > 0$.  
Therefore, the representation $\tau$ cannot be a type by Proposition \ref{keyprop}.  

Assume that $(J_0, \lambda_0)$ is attached to the simple stratum $[\Afra, n, 0, \beta_0]$ with $n>0$.  

At first, we will consider the case when $\beta \in F$.  
We have $J_0=\U(\Afra)$ by Lemma \ref{lv0twist} (i).  
Then $\rho = \lambda_0$.  
There exists a simple type $(\U(\Afra), \lambda_1)$ of level 0 and a character $\mu$ of $G$ such that $\lambda_0=\mu \otimes \lambda_1$ by Lemma \ref{lv0twist} (iii).  
We put $\tau_1=\tau \otimes \mu^{-1}$.  
Since $\tau_1$ is a type if and only if $\tau$ is a type, it is enough to show $\tau_1$ cannot be a type.  
We have $\left< \tau_1, \Ind_{K\cap{^g}\U(\Afra)}^{K} {^{gh}}\lambda_1 \right> _{K} = \left< \tau, \Ind_{K \cap {^g}\U(\Afra))}^{K} {^{gh}} \lambda_0 \right>_{K} \neq 0$.  
Therefore, by the proposition for simple types of level 0, $\tau_1$ cannot be a type.  

Then we may assume $F[\beta_0]/F$ is nontrivial.  
Since $(J_0, \lambda_0)$ is attached to the simple stratum $[\Afra, n, 0, \beta_0]$, the type $(J, \lambda)$ is a simple type attached to the simple stratum $[\Afra, n, 0, \beta] = [\Afra, n, 0, h'\beta_{0}{h'}^{-1}]$.  
Therefore, there exist representations $\kappa, \sigma$ as in Definition \ref{def_of_simple_type}.  
By our assumption, $E=F[\beta] \cong F[\beta_0]$ is an unramified extension of $F$.  
By the isomorphism $J/J^1 \cong \U(\Bfra)/\U^1(\Bfra)$, we can regard representations of $\U(\Bfra)$, trivial on $\U^1(\Bfra)$, as representations of $J$, trivial on $J^1$.  

Let $\xi$ be an irreducible summand of $\lambda|_{J \cap {^{g^{-1}}K}}$.  
Then there exists an irreducible summand $\zeta$ of $\sigma|_{J \cap {^{g^{-1}}K}}$ such that $\left< \xi, \kappa \otimes \zeta \right>_{J \cap {^{g^{-1}}K}} \neq 0$.  
Let $\mathcal{K}_0(g)=(J \cap {^{g^{-1}}}K)J^1$, and then $\mathcal{K}_0(g) \subset \mathcal{K}(g)$ and $\zeta$ can be extended to a $\mathcal{K}_0(g)$-representation which is trivial on $J^1$, since $\zeta \subset \sigma$ and $\sigma$ is trivial on $J^1$.  
The group $\mathcal{K}_0(g)$ contains $J^1$ and $\zeta$ is trivial on $J^1$, so $\zeta$ can be regarded as a $\mathcal{K}_0(g) \cap \U(\Bfra)$-representation which is trivial on $\U^1(\Bfra)$.  
Because $E/F$ is unramified and $\mathcal{K}_0(g) \subset \mathcal{K}(g)$, the image of $\mathcal{K}_0(g) \cap \U(\Bfra)$ in $\U(\Bfra)/\U^1(\Bfra)$ is sufficiently small by Proposition \ref{exam_suff_sm}.  
Since the field extension $E/F$ is nontrivial and unramified, we have $q_{E} \neq 2$.  
Then by Corollary \ref{key_to_PropertyA}, there exists an irreducible representation $\sigma'$ of $\U(\Bfra)$, trivial on $\U^1(\Bfra)$, such that $\left< \zeta , \sigma' \right>_{\U(\Bfra) \cap \mathcal{K}'(g)} \neq 0$ and $I_{B^{\times}}\left( \sigma, \sigma' \right)=\emptyset$.  
Let $\lambda' = \kappa \otimes \sigma'$.  

In the same way as in the proof of Proposition \ref{exam_intertwining}, it follows that $\lambda'$ is irreducible.  
Since $J^1 \subset \mathcal{K}_0(g) \subset J$, we have $\mathcal{K}_0(g)/J^1 \cong \left( \U(\Bfra) \cap \mathcal{K}_0(g) \right)/\U^1(\Bfra)$.  
Hence, we have $\left< \zeta , \sigma' \right>_{\mathcal{K}_0(g)}=\left< \zeta , \sigma' \right>_{\U(\Bfra) \cap \mathcal{K}_0(g)} \neq 0$.  
Therefore $\left< \zeta , \sigma' \right>_{J \cap {^{g^{-1}}}K} \geq \left< \zeta , \sigma' \right>_{\mathcal{K}_0(g)} > 0$ and $\left< \xi, \lambda' \right>_{J \cap {^{g^{-1}}}K} \geq \left< \xi, \kappa \otimes \zeta \right>_{J \cap {^{g^{-1}}K}} > 0$.  

On the other hand, by Proposition \ref{exam_intertwining}, we have $I_{G}(\lambda, \lambda') = \emptyset$.  
Therefore, by Proposition \ref{keyprop}, $\tau$ cannot be a type.  
\end{prf}

\section{Proof of main theorem}
\label{mainthm}

\begin{thm}
Let $\pi$ be an irreducible supercuspidal representation of $G$.  
Suppose $\pi$ contains a simple type $(J, \lambda)$ satisfying one of the following:  
\begin{itemize}
\item $(J, \lambda)$ is level 0,  
\item $(J, \lambda)$ is level $>0$ and attached to a simple stratum $[\Afra, n, 0, \beta]$ such that $F[\beta]$ is an unramified extension of $F$.  
\end{itemize}
Then, there exists a unique $[G,\pi]_{G}$-archetype.  

Moreover, if $K$ is a maximal compact subgroup of $G$ and $(K, \tau)$ is a $[G, \pi]_G$-type, then there exists a simple type $(J', \lambda')$ such that $J' \subset K$ and $\tau \cong \Ind_{J}^{K} \lambda'$.  
\end{thm}

\begin{prf}
Let $K$ be a maximal compact subgroup of $G$ and $(K, \tau)$ be a $[G, \pi]_{G}$-type.  
All maximal compact subgroups of $G$ are $G$-conjugate, so there exists some $x \in G$ such that ${^x}K = \GL_{m}(\ofra_{D})$.  
Since $(K, \tau)$ and $({^x}K, {^x}\tau)$ are in the same $[G, \pi]_{G}$-archetype, we may assume $K=\GL_{m}(\ofra_{D})$ and $\tau$ is defined over $K$.  
On the other hand, if $(J, \lambda)$ is attached to a simple stratum $[\Afra, n, 0, \beta]$ $({^x}J, {^x}\lambda)$ is a simple type attached to the simple stratum $[^{x}\Afra, n, 0, x\beta x^{-1}]$ for any $x \in G$, and $F[x\beta x^{-1}] \cong F[\beta]$ is an unramified extension of $F$.  
Then we may assume $\U(\Afra)$ is contained in $K$.  
By the definition of maximal simple types of level 0 and Corollary \ref{unramax}, $\Afra$ is a maximal $\ofra_F$-order in $A$ and $K=\U(\Afra)$.  
Let $\rho = \Ind_{J}^{K}\lambda$.  
Then $\rho$ is irreducible and a $[G, \pi]_{G}$-type by Corollary \ref{g1type} (this also shows the existence of a $[G, \pi]_{G}$-archetype).
Since $\tau$ is contained in $\pi|_{K}$, there exists $g \in G$ and $h \in \Kfra(\Afra)$ such that $\tau \subset \Ind_{K \cap {^g}K}^{K} {^{gh}} \rho$.  

Assume $g \notin K\Kfra(\Afra)$.  
Since $\Afra$ is maximal, there exists $g' \in Kg\Kfra(\Afra)$ such that $\left(g', Kg\Kfra(\Afra) \right)$ has property A by Lemma \ref{property_for_max}.  
Because $Kg\Kfra(\Afra) = Kg'\Kfra(\Afra)$ and $h \in \Kfra(\Afra)$, we may assume $\left(g, Kg\Kfra(\Afra) \right)$ has Property A.  
Therefore, by Proposition \ref{conclude_PropertyA}, $\tau$ cannot be a type, which is a contradiction.  
Hence $g \in K\Kfra(\Afra)=\Kfra(\Afra)$ and $\tau \subset {^{gh}} \rho$.  
These representations are irreducible, so $\tau \cong {^{gh}}\rho$.  
Then the representations $\tau$ and $\rho$ are in the same $[G,\pi]_{G}$-archetype.  
Therefore every $[G, \pi]_G$-type defined over some maximal compact subgroup of $G$ is $G$-conjugate to $(K, \rho)$, whence $[G, \pi]_{G}$-archetype is unique.  

Let $(K', \tau')$ be a $[G, \pi]_G$-type with a maximal compact subgroup $K'$ of $G$.  
Then there exists $g \in G$ such that $\tau' \cong {^g}\rho$, and then $\tau' \cong \Ind_{{^x}J}^{{^g}K} {^g}\lambda$.  
Since $(J', \lambda')=({^g}J, {^g}\lambda)$ is a simple type, the proof is completed.  
\end{prf}

\begin{rem}
If $\pi$ is an irreducible depth-zero supercuspidal representation, there exists a simple type of level 0 for $\pi$ by Proposition \ref{lv0unicity}.  
Then this theorem implies the uniqueness of $[G, \pi]_G$-archetype.  
Therefore, for inner forms of $\GL_N$, our result properly contain Latham's result \cite[Theorem 6.2]{La}.  
\end{rem}

\section{An example without the unramified assumption}

In this section, we introduce a case where there exist 
$[G, \pi]_G$-types on a maximal compact open subgroup $K$, which are not $N_{G}(K)$-conjugate, for some irreducible supercuspidal representation $\pi$.  
By this case, we find that the uniqueness of archetypes does not hold in general.  
Therefore, our example in the following is a counterexample for Latham's conjecture in \cite[Conjecture 4.4]{La}.  

Let $D$ be a quaternion algebra over $F$ and $A=\M_2(D)$.  
We take an irreducible supercuspidal representation $\pi$ of $G=GL_2(D)$ such that there exists a simple type $(J, \lambda)$ for $\pi$, attached to a simple stratum $[\Afra, n, 0, \beta]$ with the following conditions:  
\begin{itemize}
\item there exists a field extension $E=F[\beta]/E_1/F$ such that $E_1/F$ is a ramified, quadratic extension and $E/E_1$ is unramified and quadratic.  
\item $E$ is $E_1$-subalgebra in $\M_2(E_1)$ with $\ofra_{E} \subset M_2(\ofra_{E_1})$, where the embedding $\M_2(E_1) \subset M_2(D)$ is induced from some embedding $E_1 \subset D$.  
\end{itemize}
Then we have $B=\Cent_{A}(E)=E$, $\Afra=\M_2(\ofra_D)$ and $\Bfra = B \cap \Afra = \ofra_{E}$.  

Since $k_D$ and $k_{E}$ are quadratic extensions of the finite field $k_F$, the field $k_E$ identifies with $k_D$.  
We fix such an identification $k_D \cong k_E$.  
Therefore, we obtain a $k_F$-algebra inclusion 
\[
k_D=k_E=\Bfra/(\Bfra \cap \Pfra) \hookrightarrow \Afra/\Pfra = \M_2(k_D).  
\]
Note that this inclusion is not necessarily a $k_D$-algebra inclusion, so $k_E$ is not necessarily contained in the subring of diagonal matrices in $\M_2(k_D)$.  
However, there exists an element $h$ in $\U(\Afra)$ such that the image of ${}^h \Bfra$ in $\Afra/\Pfra$ is contained in the ring of diagonal matrices.  
Indeed, the above inclusion factors 
\[
\M_2(k_{F})=\M_2(k_{E_1})=\M_2(\ofra_{E_1})/\pfra_{E_1}\M_2(\ofra_{E_1}) \hookrightarrow \Afra/\Pfra = \M_2(k_D), 
\]
that is, we have a $k_F$-algebra inclusion $k_D \hookrightarrow \M_2(k_{F})$.  
Then, if $\alpha$ generates $k_E$ over $k_F$, the characteristic polynomial $\chi$ of $\alpha \in \M_2(k_{D})$ is equal to the characteristic polynomial of $\alpha \in \M_2(k_{F})$, which is the minimal polynomial of $\alpha$ over $k_F$.  
Since the polynomial $\chi$ has two different roots in $k_D$, $\alpha$ is diagonalizable as an element in $\M_2(k_D)$.  
Therefore there exists an element $x$ in $\GL_2(k_D) \cong \U(\Afra)/\U^1(\Afra)$ such that ${}^x k_{E}$ is contained the ring of diagonal matrices and a lift $h$ of $x$ in $\U(\Afra)$ satisfies the desired condition.  

We put $K=GL_2(\ofra_D)$ and $g_0=a_R(1,0)$ as in \S\ref{nice_repre}.  
Then $g_0 \notin \Kfra(\Afra)=K \Kfra(\Afra)$, that is, $Kg_0\Kfra(\Afra)$ is a non-trivial double coset of $K \backslash G/ \Kfra(\Afra)$.  
We will show ${}^{g_0h}J \subset K$.  
First, ${}^{{g_0}^{-1}}K$ contains the subgroup $K'$ in $\GL_2(\ofra_D)$ of diagonal matrices modulo $\pfra_D$.  
Since $J=\U(\Bfra)J^1 \subset \U(\Bfra) \U^1(\Afra)$, it suffices to show that ${}^h \left( \U(\Bfra) \U^1(\Afra) \right) \subset K'$.  
These groups $K'$ and ${}^h \left( \U(\Bfra) \U^1(\Afra) \right)={}^h \U(\Bfra) \U^1(\Afra)$ contain $U^1(\Afra)$, so it is enough to show that the image of ${}^h \U(\Bfra)$ in $\U(\Afra)/\U^1(\Afra)$ is contained in the counterpart of $K'$, which is already proved in the above discussion.  

We put $g=g_0 h$.  
Then $J$ and ${}^g J$ are subgroups in $K$ and we can consider $[G, \pi]_{G}$-types $\Ind_{J}^{K} \lambda$ and $\Ind_{{}^g J}^K {}^g \lambda$.  

For every element $k$ in $N_{G}(K)$, the representation ${}^k \left( \Ind_{{}^g J}^K {}^g \lambda \right)$ is not isomorphic to $\Ind_{J}^{K}\lambda$.  
Indeed, if such $k$ exists, we have $\Hom_{K}(\Ind_{J}^{K}\lambda, \Ind_{{}^{kg}J}^{K} {}^{kg}\lambda) \neq 0$ and $k'kg \in I_{G}(\lambda, \lambda)=\tilde{J}$ for some $k' \in K$.  
Because $\tilde{J} \subset \Kfra(\Afra)$ by Lemma \ref{tildeJinKA}, we obtain $g \in K\Kfra(\Afra)$.  
On the other hand, $Kg\Kfra(\Afra)=Kg_0\Kfra(\Afra)$ is a non-trivial double coset, which leads to a contradiction.

\bigbreak\bigbreak
\noindent Yuki Yamamoto\par
\noindent Graduate School of Mathematical Sciences, 
The University of Tokyo, 3--8--1 Komaba, Meguro-ku, 
Tokyo, 153--8914, Japan\par
\noindent E-mail address: \texttt{yukiymmt@ms.u-tokyo.ac.jp}


\begin{thebibliography}{99}
\bibitem{BK1} C. J. Bushnell and P. C. Kutzko, \emph{The admissible dual of $\GL(N)$ via compact open subgroups}, Annals of Mathematics Studies 129 (Princeton University Press, 1993).  
\bibitem{BK2} C. J. Bushnell and P. C. Kutzko, \emph{Smooth representations of reductive $p$-adic groups: structure theory via types}, Proc. London. Math. Soc. (3) \textbf{50} (1998), 582-634.  
\bibitem{BM} C. Breuil and A. M\'ezard, \emph{Multiplicit\'es modulaires et repr\'esentations de $\GL_{2}(\Z_p)$ et de $\Gal(\overline{\Q}_p/\Q_p)$ en $l = p$}, Duke Math. J. \textbf{115} (2002), no. 2, 205-310, With an appendix by Guy Henniart.
\bibitem{Bro} P. Broussous, \emph{Hereditary orders and embeddings of local fields in simple algebras}, J. Algebra (1) \textbf{204} (1998), 324-336.  
\bibitem{Do} A. Dotto, \emph{Breuil-M\'ezard conjectures for central division algebras}, https://arxiv.org/pdf/1808.06851, 2018.  
\bibitem{Gel} S. I. Gelfand, \emph{Representations of the general linear group over a finite field}, Lie groups and their representations (ed. I. M. Gelfand, Hilger, London, 1975) 119-132.  
\bibitem{GSZ} M. Grabitz, A. J. Silberger and E.-W. Zink, \emph{Level zero types and Hecke algebras for local central simple algebras}, J. Number Theory \textbf{77} (1) (1999), 1-26.  
\bibitem{La} P. Latham, \emph{The unicity of types for depth zero supercuspidal representations}, Represent. Theory \textbf{21} (2017), 590-610.  
\bibitem{LN} P. Latham, M. Nevins, \emph{On the unicity of types for tame toral supercuspidal representations}, https://arxiv.org/pdf/1801.06721, 2018.  
\bibitem{Pas} V. Paskunas, \emph{Unicity of types for supercuspidal representations of $\GL_{N}$}, Proc. London Math. 
Soc. (3), \textbf{91} (2005), 623-654.  
\bibitem{S1} V. S\'echerre, \emph{Repr\'esentations lisses de $\GL(m, D)$. I. caract\`eres simples}, Bull. Soc. math. France \textbf{132} (2004), no. 3, 327-396. 
\bibitem{S2} V. S\'echerre, \emph{Repr\'esentations lisses de $\GL(m, D)$. II. $\beta$-extensions}, Compos. Math. \textbf{141} (2005), no. 6, 1531-1550.  
\bibitem{S3} V. S\'echerre, \emph{Repr\'esentations lisses de $\GL(m, D)$. III. Types simples}, Ann. Sci. \'Ecole Norm. Sup. (4) \textbf{38} (2005), no. 6, 951-977.  
\bibitem{SS} V. S\'echerre and S. Stevens, \emph{Repr\'esentations lisses de $\GL(m, D)$. IV. Repr\'esentations supercuspidales}, J. Inst. Math. Jussieu \textbf{7} (2008), no. 3, 527-574.  
\bibitem{Zsig} K. Zsigmondy, \emph{Zur theorie der potenzreste}, Monatsh. Math. \textbf{3} (1892), 265-284.  
\end{thebibliography}
\end{document}